\documentclass[11pt]{article}
\usepackage{graphicx}           
\usepackage{epstopdf}
\usepackage{color}              
\usepackage{xcolor,graphicx}
\usepackage{geometry}
\usepackage{float}
\usepackage{indentfirst}        
\usepackage{amsmath,amssymb,bm} 
\usepackage{amsthm}
\usepackage{cases}              
\usepackage{setspace}
\usepackage{hyperref}
\usepackage{titlesec}
\usepackage[all]{xy}            
\usepackage{tikz}
\usetikzlibrary{arrows.meta}
\usepackage{multirow}
\usepgflibrary{arrows}
\hypersetup{dvips,ps2pdf,CJKbookmarks,hypertex,
	colorlinks=true,
	pdfpagemode=FullScreen,
	bookmarksopen=true,
	bookmarksnumbered=true}     
\usepackage{graphicx} 
\usepackage{epstopdf}
\usepackage{caption}
\usepackage{diagbox}
\usepackage{mathrsfs}
\usepackage{tikz}
\usepackage{dsfont}
\usepackage{graphics}
\usepackage{amsmath,amssymb,amsfonts,amsthm,bm,mathrsfs}
\usepackage{epsfig,graphicx,epstopdf,picinpar,subfigure,pstricks}
\usepackage{hyperref,cleveref,cite}
\usepackage{titlesec}
\usepackage{tcolorbox}
\usepackage{fancyhdr,fancyvrb}
\usepackage{xcolor,graphicx,geometry,float,indentfirst,cases,setspace,titlesec}
\usepackage{diagbox}
\usepackage{amsfonts}

\newtheorem{theorem}{Theorem}[section]

\newtheorem{lemma}[theorem]{Lemma}

\newtheorem{corollary}[theorem]{Corollary}
\newtheorem{conjecture}[theorem]{Conjecture}
\theoremstyle{remark}

\numberwithin{equation}{section}

\hfuzz=\maxdimen
\newdimen\Squaresize \Squaresize=11pt
\newdimen\Thickness \Thickness=0.7pt
\def\Square#1{\hbox{\vrule width \Thickness
   \vbox to \Squaresize{\hrule height \Thickness\vss
    \hbox to \Squaresize{\hss#1\hss}
   \vss\hrule height\Thickness}
\unskip\vrule width \Thickness} \kern-\Thickness}

\def\Vsquare#1{\vbox{\Square{$#1$}}\kern-\Thickness}

\def\moins{\raise 1pt\hbox{{$\scriptstyle -$}}}
\allowdisplaybreaks[4]

\begin{document}

\begin{center}
{\Large \bf  Inequalities for the Broken $k$-Diamond Partition Function}
\end{center}

\begin{center}
Dennis X.Q.Jia\\[6pt]
Center for Applied Mathematics\\
Tianjin University\\[6pt]
Tianjin 300072, P. R. China\\[6pt]
Email:  {jxqmail@tju.edu.cn}
\end{center}

\vspace{0.3cm} \noindent{\bf Abstract.}
In 2007, Andrews and Paule introduced the broken $k$-diamond partition function $\Delta_k(n)$, which has received a lot of researches on the arithmetic propertises. In this paper,  we prove that $D^3\log \Delta_1(n-1)>0$ for $n\geq 5$ and $D^3 \log \Delta_2(n-1)>0$ for $n\geq 7$, where $D$ is the difference operator with respect to $n$. We also conjecture that for any $k\geq 1$ and $r\geq 1$, there exists a positive integer $n_k(r)$ such that for $n\geq n_{k}(r)$, $(-1)^r D^r \log \Delta_k(n)>0$. This is analogous to the positivity of finite differences of the logarithm of the partition function, which has been proved by Chen, Wang and Xie.
    Furthermore, we obtain that both $\{\Delta_1(n)\}_{n\geq 0}$ and $\{\Delta_2(n)\}_{n\geq 0}$ satisfy the higher order Tur\'an inequalities for $n \geq 6$.

\noindent {\bf Keywords:} the higher order Tur\'{a}n inequalities, broken $k$-diamond partition function, the Jensen polynomials, finite difference
\\

\noindent {\bf AMS Classification:} 05A20, 11P82
\\


\section{INTRODUCTION}
\begin{figure}[htbp]
	\begin{center}
		\begin{tikzpicture}[>=Stealth,scale=0.8]	
		\coordinate[label=above:$a_2$] (a2) at (0,1.4);
		\draw(a2) circle(.1);
		
		\coordinate[label=below:$a_3$] (a3) at (0,-1.4);
		\draw (a3) circle(.1);
			
		\coordinate[label=left:$a_1$] (a1) at (-1.4,0);
		\draw (a1) circle(.1);
		
		\coordinate [label=right:$a_4$] (a4) at (1.4,0);
		\draw (a4) circle(.1);
		
		\draw[line width=0.5pt](-1.33,0.07)--(-0.07,1.33);
		\draw[-latex, line width=0.5pt] (-1.33,0.07)--(-0.6,0.8);
		
		\draw[line width=0.5pt](-1.33,-0.07)--(-0.07,-1.33);
		\draw[-latex, line width=0.5pt] (-1.33,-0.07)--(-0.6,-0.8);

		\draw[line width=0.5pt](0.07,1.33)--(1.33,0.07);
		\draw[-latex, line width=0.5pt](0.07,1.33)--(0.8,0.6);
		
		\draw[line width=0.5pt](0.07,-1.33)--(1.33,-0.07);
		\draw[-latex, line width=0.5pt](0.07,-1.33)--(0.8,-0.6);
		
\end{tikzpicture}
		\caption{A plane partition satisfied \eqref{eq-part}.}
		\label{fig-plane-partition}
	\end{center}
\end{figure}
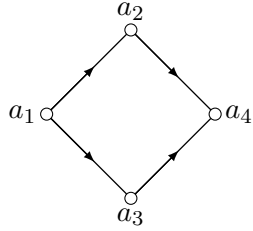

\begin{figure}[h]
	\begin{center}
		\begin{tikzpicture}[>=Stealth,scale=0.8]	
		
		\coordinate[label=above:$a_2$] (a2) at (0,1.4);
		\draw(a2) circle(.1);
		
		\coordinate[label=below:$a_3$] (a3) at (0,-1.4);
		\draw (a3) circle(.1);
			
		\coordinate[label=left:$a_1$] (a1) at (-1.4,0);
		\draw (a1) circle(.1);
		
		\coordinate (a4) at (1.4,0);
		\draw (a4) circle(.1);
		\node at (1.4, 0.1) [above]{$a_4$};
		
		\draw[line width=0.5pt](-1.33,0.07)--(-0.07,1.33);
		\draw[-latex, line width=0.5pt] (-1.33,0.07)--(-0.6,0.8);
		
		\draw[line width=0.5pt](-1.33,-0.07)--(-0.07,-1.33);
		\draw[-latex, line width=0.5pt] (-1.33,-0.07)--(-0.6,-0.8);

		\draw[line width=0.5pt](0.07,1.33)--(1.33,0.07);
		\draw[-latex, line width=0.5pt](0.07,1.33)--(0.8,0.6);
		
		\draw[line width=0.5pt](0.07,-1.33)--(1.33,-0.07);
		\draw[-latex, line width=0.5pt](0.07,-1.33)--(0.8,-0.6);

		\coordinate[label=above:$a_5$] (a5) at (2.8,1.4);
		\draw(a5) circle(.1);
		
		\coordinate[label=below:$a_6$] (a6) at (2.8,-1.4);
		\draw (a6) circle(.1);			

		\coordinate  (a7) at (4.2,0);
		\draw (a7) circle(.1);
		\node at (4.2,0.1)[above]{$a_7$};
		
		\draw[line width=0.5pt](1.47,0.07)--(2.73,1.33);
		\draw[-latex, line width=0.5pt] (1.47,0.07)--(2.2,0.8);
		
		\draw[line width=0.5pt](1.47,-0.07)--(2.73,-1.33);
		\draw[-latex, line width=0.5pt] (1.47,-0.07)--(2.2,-0.8);

		\draw[line width=0.5pt](2.87,1.33)--(4.13,0.07);
		\draw[-latex, line width=0.5pt](2.87,1.33)--(3.6,0.6);
		
		\draw[line width=0.5pt](2.87,-1.33)--(4.13,-0.07);
		\draw[-latex, line width=0.5pt](2.87,-1.33)--(3.6,-0.6);
		
		\coordinate[label=above:$a_8$] (a8) at (5.6,1.4);
		\draw(a8) circle(.1);
		
		\coordinate[label=below:$a_9$] (a9) at (5.6,-1.4);
		\draw (a9) circle(.1);			

		\coordinate  (a10) at (7,0);
		\draw (a10) circle(.1);
		\node at (7,0.1)[above]{$a_{10}$};
		
		\draw[line width=0.5pt](4.27,0.07)--(5.53,1.33);
		\draw[-latex, line width=0.5pt] (4.27,0.07)--(5,0.8);
		
		\draw[line width=0.5pt](4.27,-0.07)--(5.53,-1.33);
		\draw[-latex, line width=0.5pt] (4.27,-0.07)--(5,-0.8);

		\draw[line width=0.5pt](5.67,1.33)--(6.93,0.07);
		\draw[-latex, line width=0.5pt](5.67,1.33)--(6.4,0.6);
		
		\draw[line width=0.5pt](5.67,-1.33)--(6.93,-0.07);
		\draw[-latex, line width=0.5pt](5.67,-1.33)--(6.4,-0.6);

		\coordinate[label=above:$a_{3\ell-1}$] (a3n-1) at (10.5,1.4);
		\draw(a3n-1) circle(.1);
		
		\coordinate[label=below:$a_{3\ell}$] (a3n) at (10.5,-1.4);
		\draw (a3n) circle(.1);

		\coordinate (a3n-2) at (9.1,0);
		\draw (a3n-2) circle(.1);
		\node at (8.7,0.1) [above]{$a_{3\ell-2}$};
		
		\node at (7.3,0)[]{$\cdot$};
		\node at (7.6,0)[]{$\cdot$};
		\node at (7.9,0)[]{$\cdot$};
		\node at (8.2,0)[]{$\cdot$};
		\node at (8.5,0)[]{$\cdot$};
		\node at (8.8,0)[]{$\cdot$};
		
		\coordinate[label=right:$a_{3\ell+1}$] (a3n+1) at (11.9,0);
		\draw (a3n+1) circle(.1);

		\draw[line width=0.5pt](9.17,0.07)--(10.43,1.33);
		\draw[-latex, line width=0.5pt] (9.17,0.07)--(9.9,0.8);
		
		\draw[line width=0.5pt](9.17,-0.07)--(10.43,-1.33);
		\draw[-latex, line width=0.5pt] (9.17,-0.07)--(9.9,-0.8);

		\draw[line width=0.5pt](10.57,1.33)--(11.83,0.07);
		\draw[-latex, line width=0.5pt](10.57,1.33)--(11.3,0.6);
		
		\draw[line width=0.5pt](10.57,-1.33)--(11.83,-0.07);
		\draw[-latex, line width=0.5pt](10.57,-1.33)--(11.3,-0.6);
	\end{tikzpicture}
		\caption{A plane partition diamond of length $\ell$.}
		\label{fig-plane-partition-diamonds}
	\end{center}
\end{figure}
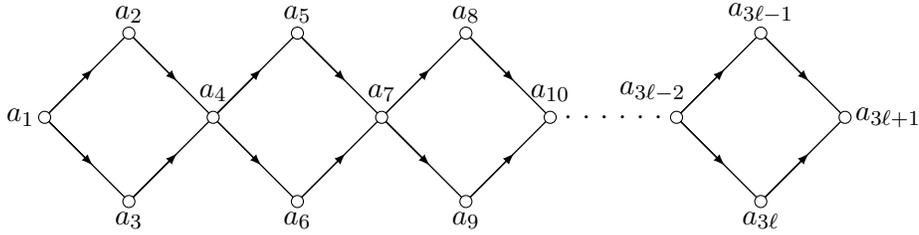

The classical plane partition, introduced by MacMahon in his famous book  "Combinatory Analysis'' \cite{MacMahon-1915}, is the situation where the non-negative integer parts $a_i$ of the partition are placed at the vertices of a square such that the following relations are satisfied:
\begin{equation}\label{eq-part}
	a_1\geq a_2, ~a_1\geq a_3, ~a_2\geq a_4, ~\text{and} ~ a_3\geq a_4.
\end{equation}
With the arrows describing the relation ``$\geq$'', one can represent these partitions as directed graphs. For example, Figure \ref{fig-plane-partition} represents the relation \eqref{eq-part}.
Here and throughout the  following an arrow pointing from $a_i$ to $a_j$ is always interpreted as $a_i \geq a_j$. 
Note that the arrows can also be interpreted as a ``ﬂow'', so one can regard $a_1$ as the ``source'' and $a_4$ as the ``sink''.
A plane partition diamond of length $\ell$, named by
Andrews, Paule and Riese \cite{Andrews-Paule-Riese-2001}, is shown in Figure \ref{fig-plane-partition-diamonds}.
Later, Andrews and Paule \cite{Andrews-Paule-2007} introduced the idea of $k$-elongated partition diamonds.
 A $k$-elongated diamond of length $1$ and length $\ell$ are depicted in Figure \ref{fig-k-elo-1} and Figure \ref{fig-k-elo-n}, respectively. 
\begin{figure}[h]
	\begin{center}
		\begin{tikzpicture}[>=Stealth,scale=0.8]	
			\draw (0,0) circle(.1);
			\node at (0,0) [above] {$a_3$};	
			\draw (1.4,0) circle(.1);
			\node at (1.4,0) [above] {$a_5$};	
			\draw (2.8,0) circle(.1);
			\node at (2.8,0) [above] {$a_7$};
			
			\node at (3.1,0) []{$\cdot$};
			\node at (3.4,0) []{$\cdot$};
			\node at (3.7,0) []{$\cdot$};
			\node at (4,0) []{$\cdot$};
			\node at (4.3,0) []{$\cdot$};
			\node at (4.6,0) []{$\cdot$};

			\draw (4.9,0) circle(.1);
			\node at (4.9,0) [above] {$a_{2k-1}$};
			\draw (6.3,0) circle(.1);
			\node at (6.3,0) [above] {$a_{2k+1}$};	
			\draw (7,-1.4) circle(.1);
			\node at (7,-1.4) [right] {$a_{2k+2}$};	
			
			\draw (-0.7,-1.4) circle(.1);
			\node at (-0.7,-1.4) [left] {$a_1$};	
			\draw (0,-2.8) circle(.1);
			\node at (0,-2.8) [below] {$a_2$};
			\draw (1.4,-2.8) circle(.1);
			\node at (1.4,-2.8) [below] {$a_4$};		
			\draw (2.8,-2.8) circle(.1);
			\node at (2.8,-2.8) [below] {$a_6$};
			
			\node at (3.1,-2.8) []{$\cdot$};
			\node at (3.4,-2.8) []{$\cdot$};
			\node at (3.7,-2.8) []{$\cdot$};
			\node at (4,-2.8) []{$\cdot$};
			\node at (4.3,-2.8) []{$\cdot$};
			\node at (4.6,-2.8) []{$\cdot$};
			
			\draw (4.9,-2.8) circle(.1);
			\node at (4.9,-2.8) [below] {$a_{2k-2}$};
			\draw (6.3,-2.8) circle(.1);
			\node at (6.3,-2.8) [below] {$a_{2k}$};
			
			\draw[-latex,line width=0.6pt](0.1,0)--(0.8,0);
			\draw[line width=0.6pt] (0.7,0)--(1.3,0);
			\draw[-latex,line width=0.6pt](1.5,0)--(2.2,0);
			\draw[line width=0.6pt] (2.1,0)--(2.7,0);
			\draw[-latex,line width=0.6pt](5,0)--(5.7,0);
			\draw[line width=0.6pt] (5.6,0)--(6.2,0);
			
			\draw[-latex,line width=0.6pt](6.345,-0.09)--(6.7,-0.8);
			\draw[line width=0.6pt] (6.6,-0.6)--(6.955,-1.31);
			
			\draw[-latex,line width=0.6pt](-0.655,-1.31)--(-0.3,-0.6);
			\draw[line width=0.6pt] (-0.4,-0.8)--(-0.045,-0.09);
			
			\draw[-latex,line width=0.6pt](-0.655,-1.49)--(-0.3,-2.2);
			\draw[line width=0.6pt] (-0.4,-2)--(-0.045,-2.71);
			
			\draw[-latex,line width=0.6pt](0.1,-2.8)--(0.8,-2.8);
			\draw[line width=0.6pt] (0.7,-2.8)--(1.3,-2.8);
			
			\draw[-latex,line width=0.6pt](1.5,-2.8)--(2.2,-2.8);
			\draw[line width=0.6pt] (2.1,-2.8)--(2.7,-2.8);
			
			\draw[-latex,line width=0.6pt](5,-2.8)--(5.7,-2.8);
			\draw[line width=0.6pt] (5.6,-2.8)--(6.2,-2.8);
			
			\draw[-latex,line width=0.6pt](6.345,-2.71)--(6.7,-2);
			\draw[line width=0.6pt] (6.6,-2.2)--(6.955,-1.49);
			
			\draw[-latex,line width=0.6pt](0.045,-0.09)--(0.4,-0.8);
			\draw[line width=0.6pt] (0.3,-0.6)--(1.355,-2.71);
			
			\draw[-latex,line width=0.6pt](1.445,-0.09)--(1.8,-0.8);
			\draw[line width=0.6pt] (1.7,-0.6)--(2.755,-2.71);
			
			\draw[-latex,line width=0.6pt](4.945,-0.09)--(5.3,-0.8);
			\draw[line width=0.6pt] (5.2,-0.6)--(6.255,-2.71);

			\draw[-latex,line width=0.6pt](0.045,-2.71)--(0.4,-2);
			\draw[line width=0.6pt] (0.3,-2.2)--(1.355,-0.09);
			
			\draw[-latex,line width=0.6pt](1.445,-2.71)--(1.8,-2);
			\draw[line width=0.6pt] (1.7,-2.2)--(2.755,-0.09);
			
			\draw[-latex,line width=0.6pt](4.945,-2.71)--(5.3,-2);
			\draw[line width=0.6pt] (5.2,-2.2)--(6.255,-0.09);
			
		\end{tikzpicture}
		
		\caption{ A $k$-elongated partition diamond of length $1$.}
		\label{fig-k-elo-1}
	\end{center}
\end{figure}
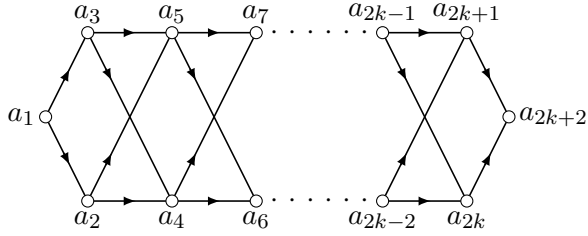
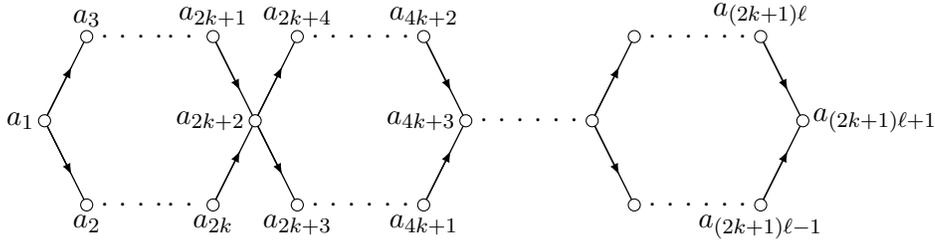
\begin{figure}[h]
	\begin{center}
		\begin{tikzpicture}[>=Stealth,scale=0.8]	
		
			\coordinate[label=left:$a_1$] (a1) at (0,0);
			\draw (a1) circle(.1);
			
			\coordinate[label=above:$a_3$] (a3) at (0.7,1.4);
			\draw (a3) circle(.1);
			\coordinate [label=below:$a_2$] (a2) at (0.7,-1.4);
			\draw (a2) circle(.1);
			\node at (1,1.4) {$\cdot$};	
			\node at (1.3,1.4) {$\cdot$};	
			\node at (1.6,1.4) {$\cdot$};	
			\node at (1.9,1.4) {$\cdot$};	
			\node at (2.1,1.4) {$\cdot$};	
			\node at (2.4,1.4) {$\cdot$};
				
			\coordinate [label=above:$a_{2k+1}$] (a2k+1) at (2.8,1.4);
			\draw (a2k+1) circle(.1);
			\coordinate [label=below:$a_{2k}$] (a2k) at (2.8,-1.4);
			\draw (a2k) circle(.1);
			\node at (1,-1.4) {$\cdot$};	
			\node at (1.3,-1.4) {$\cdot$};	
			\node at (1.6,-1.4) {$\cdot$};	
			\node at (1.9,-1.4) {$\cdot$};	
			\node at (2.1,-1.4) {$\cdot$};	
			\node at (2.4,-1.4) {$\cdot$};
			
			\coordinate [label=left:$a_{2k+2}$] (a2k+2) at (3.5,0);
			\draw (a2k+2) circle(.1);

			\draw[-latex,line width=0.5pt] (0.045,0.09)--(0.45,0.9);
			\draw[line width=0.5pt] (0.045,0.09)--(0.655,1.31);
			
			\draw[-latex,line width=0.5pt] (0.045,-0.09)--(0.45,-0.9);
			\draw[line width=0.5pt] (0.045,-0.09)--(0.655,-1.31);
			
			\draw[line width=0.5pt] (2.845,1.31)--(3.455,0.09);
			\draw[-latex,line width=0.5pt] (2.845,1.31)--(3.25,0.5);
			
			\draw[-latex,line width=0.5pt] (2.845,-1.31)--(3.25,-0.5);
			\draw[line width=0.5pt] (2.845,-1.31)--(3.455,-0.09);

			\coordinate[label=above:$a_{2k+4}$] (a2k+4) at (4.2,1.4);
			\draw (a2k+4) circle(.1);
			\coordinate [label=below:$a_{2k+3}$] (a2k+3) at (4.2,-1.4);
			\draw (a2k+3) circle(.1);
			\node at (4.5,1.4) {$\cdot$};	
			\node at (4.8,1.4) {$\cdot$};	
			\node at (5.1,1.4) {$\cdot$};	
			\node at (5.4,1.4) {$\cdot$};	
			\node at (5.7,1.4) {$\cdot$};	
			\node at (6,1.4) {$\cdot$};
				
			\coordinate [label=above:$a_{4k+2}$] (a4k+2) at (6.3,1.4);
			\draw (a4k+2) circle(.1);
			\coordinate [label=below:$a_{4k+1}$] (a4k+1) at (6.3,-1.4);
			\draw (a4k+1) circle(.1);
			\node at (4.5,-1.4) {$\cdot$};	
			\node at (4.8,-1.4) {$\cdot$};	
			\node at (5.1,-1.4) {$\cdot$};	
			\node at (5.4,-1.4) {$\cdot$};	
			\node at (5.7,-1.4) {$\cdot$};	
			\node at (6,-1.4) {$\cdot$};
			
			\coordinate [label=left:$a_{4k+3}$] (a4k+3) at (7,0);
			\draw (a4k+3) circle(.1);

			\draw[-latex,line width=0.5pt] (3.545,0.09)--(3.95,0.9);
			\draw[line width=0.5pt] (3.545,0.09)--(4.155,1.31);
			
			\draw[-latex,line width=0.5pt] (3.545,-0.09)--(3.95,-0.9);
			\draw[line width=0.5pt] (3.545,-0.09)--(4.155,-1.31);
			
			\draw[line width=0.5pt] (6.345,1.31)--(6.955,0.09);
			\draw[-latex,line width=0.5pt] (6.345,1.31)--(6.75,0.5);
			
			\draw[-latex,line width=0.5pt] (6.345,-1.31)--(6.75,-0.5);
			\draw[line width=0.5pt] (6.345,-1.31)--(6.955,-0.09);
			
			\node at (7.3,0) {$\cdot$};	
			\node at (7.6,0) {$\cdot$};	
			\node at (7.9,0) {$\cdot$};	
			\node at (8.2,0) {$\cdot$};	
			\node at (8.5,0) {$\cdot$};	
			\node at (8.8,0) {$\cdot$};

			\draw (9.1,0) circle(.1);

			\draw (9.8,1.4) circle(.1);
			
			\draw (9.8,-1.4) circle(.1);

			\node at (10.1,1.4) {$\cdot$};	
			\node at (10.4,1.4) {$\cdot$};	
			\node at (10.7,1.4) {$\cdot$};	
			\node at (11,1.4) {$\cdot$};	
			\node at (11.3,1.4) {$\cdot$};	
			\node at (11.6,1.4) {$\cdot$};
				
			\coordinate [label=above:$a_{(2k+1)\ell}$] (a2k+1n) at (11.9,1.4);
			\draw (a2k+1n) circle(.1);
			
			\coordinate [label=below:$a_{(2k+1)\ell-1}$] (a2k+1n-1) at (11.9,-1.4);
			\draw (a2k+1n-1) circle(.1);
			
			\node at (10.1,-1.4) {$\cdot$};	
			\node at (10.4,-1.4) {$\cdot$};	
			\node at (10.7,-1.4) {$\cdot$};	
			\node at (11,-1.4) {$\cdot$};	
			\node at (11.3,-1.4) {$\cdot$};	
			\node at (11.6,-1.4) {$\cdot$};
			
			\coordinate [label=right:$a_{(2k+1)\ell+1}$] (a2k+1n+1) at (12.6,0);
			\draw (a2k+1n+1) circle(.1);

			\draw[-latex,line width=0.5pt] (9.145,0.09)--(9.55,0.9);
			\draw[line width=0.5pt] (9.145,0.09)--(9.755,1.31);
			
			\draw[-latex,line width=0.5pt] (9.145,-0.09)--(9.55,-0.9);
			\draw[line width=0.5pt] (9.145,-0.09)--(9.755,-1.31);
			
			\draw[line width=0.5pt] (11.945,1.31)--(12.555,0.09);
			\draw[-latex,line width=0.5pt] (11.945,1.31)--(12.35,0.5);
			
			\draw[-latex,line width=0.5pt] (11.945,-1.31)--(12.35,-0.5);
			\draw[line width=0.5pt] (11.945,-1.31)--(12.555,-0.09);

		\end{tikzpicture}
		\caption{ A $k$-elongated partition diamond of length $\ell$.}
		\label{fig-k-elo-n}
	\end{center}
\end{figure}
 Moreover, they also introduced a new variation of plane partitions diamond, called  broken $k$-diamond partitions. 
 A broken $k$-diamond partition of length $2\ell$ is shown in Figure \ref{fig-broken-k}. 
 As shown in the figure, a broken $k$-diamond consists of two separated $k$-elongated diamond partitions of length $\ell$,
 in which one of them  the source $b_1$ is deleted.
 \begin{figure}[h]
	\begin{center}
		\begin{tikzpicture}[>=Stealth,scale=0.6]	
			
			\coordinate [label=left:$a_{1}$] (a1) at (2.1,0);
			\draw (a1) circle(.1);
			
			\coordinate[label=above:$a_{5}$] (a5) at (4.2,1.4);
			\draw (a5) circle(.1);
			
			\coordinate [label=below:$a_{4}$] (a4) at (4.2,-1.4);
			\draw (a4) circle(.1);
			
			\coordinate[label=above:$a_{3}$] (a3) at (2.8,1.4);
			\draw (a3) circle(.1);
			
			\coordinate[label=below:$a_{2}$] (a2) at (2.8,-1.4);
			\draw (a2) circle(.1);
			
			\draw[line width=0.5pt] (2.9,1.4)--(4.1,1.4);
			\draw[-latex,line width=0.5pt] (2.9,1.4)--(3.7,1.4);
			
			\draw[line width=0.5pt] (2.9,-1.4)--(4.1,-1.4);
			\draw[-latex,line width=0.5pt] (2.9,-1.4)--(3.7,-1.4);
			
			\node at (4.5,1.4) {$\cdot$};	
			\node at (4.8,1.4) {$\cdot$};	
			\node at (5.1,1.4) {$\cdot$};	
			\node at (5.4,1.4) {$\cdot$};	
			\node at (5.7,1.4) {$\cdot$};	
			\node at (6,1.4) {$\cdot$};
				
			\coordinate [label=above:$a_{2k+1}$] (a2k+1) at (6.3,1.4);
			\draw (a2k+1) circle(.1);
			
			\coordinate [label=below:$a_{2k}$] (a2k) at (6.3,-1.4);
			\draw (a2k) circle(.1);
			
			\node at (4.5,-1.4) {$\cdot$};	
			\node at (4.8,-1.4) {$\cdot$};	
			\node at (5.1,-1.4) {$\cdot$};	
			\node at (5.4,-1.4) {$\cdot$};	
			\node at (5.7,-1.4) {$\cdot$};	
			\node at (6,-1.4) {$\cdot$};
			
			\coordinate [label=left:$a_{2k+2}$] (a2k+2) at (7,0);
			\draw (a2k+2) circle(.1);

			\draw[line width=0.5pt] (2.845,-1.31)--(4.155,1.31);
			\draw[-latex,line width=0.5pt] (2.845,-1.31)--(3.3,-0.4);

			\draw[line width=0.5pt] (2.845,1.31)--(4.155,-1.31);
			\draw[-latex,line width=0.5pt] (2.845,1.31)--(3.3,0.4);
			
			\draw[-latex,line width=0.5pt] (2.145,0.09)--(2.55,0.9);
			\draw[line width=0.5pt] (2.145,0.09)--(2.755,1.31);
			
		    \draw[-latex,line width=0.5pt] (2.145,-0.09)--(2.55,-0.9);
			\draw[line width=0.5pt] (2.145,-0.09)--(2.755,-1.31);
			
			\draw[line width=0.5pt] (6.345,1.31)--(6.955,0.09);
			\draw[-latex,line width=0.5pt] (6.345,1.31)--(6.75,0.5);
			
			\draw[-latex,line width=0.5pt] (6.345,-1.31)--(6.75,-0.5);
			\draw[line width=0.5pt] (6.345,-1.31)--(6.955,-0.09);
			
			\node at (7.3,0) {$\cdot$};	
			\node at (7.6,0) {$\cdot$};	
			\node at (7.9,0) {$\cdot$};

			\draw (8.2,0) circle(.1);

			\draw (8.9,1.4) circle(.1);
			
			\draw (8.9,-1.4) circle(.1);

			\node at (9.2,1.4) {$\cdot$};	
			\node at (9.5,1.4) {$\cdot$};	
			\node at (9.8,1.4) {$\cdot$};	
			\node at (10.1,1.4) {$\cdot$};	
			\node at (10.4,1.4) {$\cdot$};	
			\node at (10.7,1.4) {$\cdot$};
				
			\coordinate [label=above:$a_{(2k+1)\ell}$] (a2k+1n) at (11,1.4);
			\draw (a2k+1n) circle(.1);
			
			\coordinate [label=below:$a_{(2k+1)\ell-1}$] (a2k+1n-1) at (11,-1.4);
			\draw (a2k+1n-1) circle(.1);
			
			\node at (9.2,-1.4) {$\cdot$};	
			\node at (9.5,-1.4) {$\cdot$};	
			\node at (9.8,-1.4) {$\cdot$};	
			\node at (10.1,-1.4) {$\cdot$};	
			\node at (10.4,-1.4) {$\cdot$};	
			\node at (10.7,-1.4) {$\cdot$};
			
			\coordinate [label=left:$a_{(2k+1)\ell+1}$] (a2k+1n+1) at (11.7,0);
			\draw (a2k+1n+1) circle(.1);

			\draw[-latex,line width=0.5pt] (8.245,0.09)--(8.65,0.9);
			\draw[line width=0.5pt] (8.245,0.09)--(8.855,1.31);
			
			\draw[-latex,line width=0.5pt] (8.245,-0.09)--(8.65,-0.9);
			\draw[line width=0.5pt] (8.245,-0.09)--(8.855,-1.31);
			
			\draw[line width=0.5pt] (11.045,1.31)--(11.655,0.09);
			\draw[-latex,line width=0.5pt] (11.045,1.31)--(11.45,0.5);
			
			\draw[-latex,line width=0.5pt] (11.045,-1.31)--(11.45,-0.5);
			\draw[line width=0.5pt] (11.045,-1.31)--(11.655,-0.09);


			
			\coordinate[label=above:$b_{5}$] (b5) at (-1.2,1.4);
			\draw (b5) circle(.1);
			
			\coordinate [label=below:$b_{4}$] (b4) at (-1.2,-1.4);
			\draw (b4) circle(.1);
			
			\coordinate[label=above:$b_{3}$] (b3) at (0.2,1.4);
			\draw (b3) circle(.1);
			
			\coordinate[label=below:$b_{2}$] (b2) at (0.2,-1.4);
			\draw (b2) circle(.1);
			
			\draw[line width=0.5pt] (0.1,1.4)--(-1.1,1.4);
			\draw[-latex,line width=0.5pt] (0.1,1.4)--(-0.7,1.4);
			
			\draw[line width=0.5pt] (0.1,-1.4)--(-1.1,-1.4);
			\draw[-latex,line width=0.5pt] (0.1,-1.4)--(-0.7,-1.4);
			
			\node at (-1.5,1.4) {$\cdot$};	
			\node at (-1.8,1.4) {$\cdot$};	
			\node at (-2.1,1.4) {$\cdot$};	
			\node at (-2.4,1.4) {$\cdot$};	
			\node at (-2.7,1.4) {$\cdot$};	
			\node at (-3,1.4) {$\cdot$};
				
			\coordinate [label=above:$b_{2k+1}$] (b2k+1) at (-3.3,1.4);
			\draw (b2k+1) circle(.1);
			
			\coordinate [label=below:$b_{2k}$] (b2k) at (-3.3,-1.4);
			\draw (b2k) circle(.1);
			
			\node at (-1.5,-1.4) {$\cdot$};	
			\node at (-1.8,-1.4) {$\cdot$};	
			\node at (-2.1,-1.4) {$\cdot$};	
			\node at (-2.4,-1.4) {$\cdot$};	
			\node at (-2.7,-1.4) {$\cdot$};	
			\node at (-3,-1.4) {$\cdot$};
			
			\coordinate [label=right:$b_{2k+2}$] (b2k+2) at (-4,0);
			\draw (b2k+2) circle(.1);

			\draw[line width=0.5pt] (0.155,-1.31)--(-1.155,1.31);
			\draw[-latex,line width=0.5pt] (0.155,-1.31)--(-0.3,-0.4);

			\draw[line width=0.5pt] (0.155,1.31)--(-1.155,-1.31);
			\draw[-latex,line width=0.5pt] (0.155,1.31)--(-0.3,0.4);

			\draw[line width=0.5pt] (-3.345,1.31)--(-3.955,0.09);
			\draw[-latex,line width=0.5pt] (-3.345,1.31)--(-3.75,0.5);
			
			\draw[-latex,line width=0.5pt] (-3.345,-1.31)--(-3.75,-0.5);
			\draw[line width=0.5pt] (-3.345,-1.31)--(-3.955,-0.09);
			
			\node at (-4.3,0) {$\cdot$};	
			\node at (-4.6,0) {$\cdot$};	
			\node at (-4.9,0) {$\cdot$};

			\draw (-5.2,0) circle(.1);

			\draw (-5.9,1.4) circle(.1);
			
			\draw (-5.9,-1.4) circle(.1);

			\node at (-6.2,1.4) {$\cdot$};	
			\node at (-6.5,1.4) {$\cdot$};	
			\node at (-6.8,1.4) {$\cdot$};	
			\node at (-7.1,1.4) {$\cdot$};	
			\node at (-7.4,1.4) {$\cdot$};	
			\node at (-7.7,1.4) {$\cdot$};
				
			\coordinate [label=above:$b_{(2k+1)\ell}$] (b2k+1n) at (-8,1.4);
			\draw (b2k+1n) circle(.1);
			
			\coordinate [label=below:$b_{(2k+1)\ell-1}$] (b2k+1n-1) at (-8,-1.4);
			\draw (b2k+1n-1) circle(.1);
			
			\node at (-6.2,-1.4) {$\cdot$};	
			\node at (-6.5,-1.4) {$\cdot$};	
			\node at (-6.8,-1.4) {$\cdot$};	
			\node at (-7.1,-1.4) {$\cdot$};	
			\node at (-7.4,-1.4) {$\cdot$};	
			\node at (-7.7,-1.4) {$\cdot$};
			
			\coordinate [label=left:$b_{(2k+1)+1}$] (b2k+1n+1) at (-8.7,0);
			\draw (b2k+1n+1) circle(.1);

			\draw[-latex,line width=0.5pt] (-5.245,0.09)--(-5.65,0.9);
			\draw[line width=0.5pt] (-5.245,0.09)--(-5.855,1.31);
			
			\draw[-latex,line width=0.5pt] (-5.245,-0.09)--(-5.65,-0.9);
			\draw[line width=0.5pt] (-5.245,-0.09)--(-5.855,-1.31);
			
			\draw[line width=0.5pt] (-8.045,1.31)--(-8.655,0.09);
			\draw[-latex,line width=0.5pt] (-8.045,1.31)--(-8.45,0.5);
			
			\draw[-latex,line width=0.5pt] (-8.045,-1.31)--(-8.45,-0.5);
			\draw[line width=0.5pt] (-8.045,-1.31)--(-8.655,-0.09);
				
		\end{tikzpicture}
		\caption{ A broken $k$-diamond of length $2\ell$.}
		\label{fig-broken-k}
	\end{center}
\end{figure}
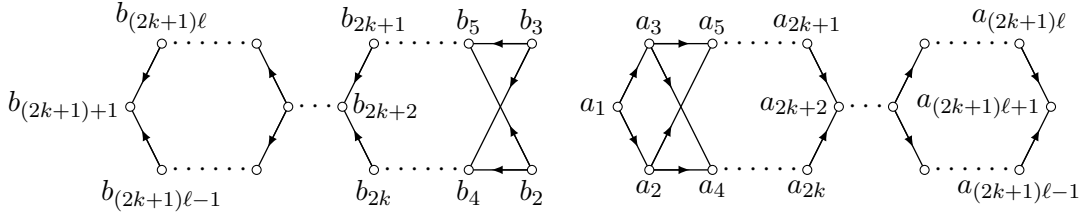
 For any integer  $k\geq 0$,  let $\Delta_k(n)$ denote the 
 the number of broken $k$-diamond partitions of a non-negative integer $n$.
  Andrews and Paule \cite{Andrews-Paule-2007} proved that 
\begin{align*}
	\sum_{n=0}^{\infty} \Delta_k(n) q^n =\prod_{n=1}^{\infty} \frac{(1-q^{2n})(1-q^{(2k+1)n})}{(1-q^n)^3(1-q^{(4k+2)n})}.
\end{align*}
Note that the generating function for $\Delta_k(n)$ is a modular form. More precisely, $\Delta_k(n)$ are the coefficients of a modular function with respect to $\Gamma_0(4k+2)$. 
The modular aspect led to some various arithmetic theorems and properties.
For instance, many Ramanujan-like congruences satisfied by $\Delta_k(n)$ has been proved by many authors, see 
\cite{Andrews-Paule-2007,Chan-2008,
Chen-Fan-Yu-2014,
Hirschhorn-2018,
Hirschhorn-Sellers-2007,
Jameson-2013,
Paule-Radu-2010,Xiong-2011} and so on.

The objective of this paper is to prove both $\{\Delta_1(n)\}_{n\geq 6}$ and $\{\Delta_2(n)\}_{n\geq 6}$ satisfy
 the higher order Tur\'an inequalities.
  The higher order Tur\'an inequalities (also called the cubic inequalities) and the Tur\'an inequalities (also named the Newton inequalities) arise in the
	study of the Maclaurin coefficients of real entire functions in the
	Laguerre-P\'olya class. We refer the interested readers to  \cite{Levin-1980} and \cite{Rahman-2002}. 
	A real sequence $\{a_n\}
	_{n\geq 0}$ is said to satisfy the Tu\'an inequalities or to be log-concave
	if for $n\geq 1$,
	\begin{align}
	a_n^2-a_{n-1}a_{n+1}\geq 0;
	\end{align}
	it is said to satisfy the higher order Tur\'an inequalities if for $n\geq 1$,
	\begin{align}
		4(a_n^2-a_{n-1}a_{n+1})(a_{n+1}^2-a_na_{n+1})-(a_na_{n+1}-a_{n-1}a_{n+2})^2\geq 0;
	\end{align}
	see \cite{Craven-1989,Csordas-1986,Dimitrov-1998,Niculescu-2000, Rosset-1989,Wagner-1977}. The Tur\'an inequalities and higher order Tur\'an inequalities are closely related with the Jensen polynomials. The Jensen polynomials of degree $d$ and shift $n$ associated to the sequence $\{a_n\}_{n\geq 0}$ are defined by 
	\begin{equation}
		J_{a}^{d,n}(X)=\sum_{j=0}^{d}\binom{d}{j}a_{n+j} X^{j}.
	\end{equation}
	More properties of the Jensen polynomials can be found in \cite{Craven-1989,Csordas-1986,Csordas-1990}.
	We note that the quadratic Jensen polynomial associated to $\{a_n\}_{n\geq 0}$ of shift $n-1$  is  
	\begin{equation}
		J_{a}^{2,n-1}(X)=a_{n-1}+2  a_{n}X+a_{n+1}X^2.
	\end{equation}
	Then one sees that $\{a_n\}_{n\geq 0}$ is log-concave at $n$ if and only
	if $J_{a}^{2,n-1}(X)$ is hyperbolic. 
	A real polynomial is said to be hyperbolic if all of its zeros are real.
	It is said in general that
	a sequence $\{a_n\}_{n\geq 0}$ satisfies the order $d$ Tur\'an
	inequality at $n$ if and only if $J_{a}^{d,n-1}(X)$ is hyperbolic.

	Recently, many combinatorial sequences  satisfying the (higher) Tur\'an inequalities have been proved  by a lot of mathematicians. 
	For example, Nicolas
	\cite{Nicolas-1978} and DeSalvo and Pak \cite{DeSalvo-Pak-2015} proved that the partition
	function $p(n)$ is log-concave for $n\geq 26$, where $p(n)$ is the
	number of partitions of $n$. Chen \cite{Chen-2017} conjectured that the partition function $p(n)$
	satisfies the higher order Tur\'an inequalities for $n\geq 95$, which
	has been proved by Chen, Jia and Wang \cite{Chen-Jia-Wang-2019} and Larson and Wagner \cite{Larson-2019}. In fact, for sufficiently large $n$ and any $d \geq 1$, Griffin, Ono, Rolen and Zagier \cite{Griffin-Ono-Rolen-Zagier-2019} proved the partition function satisfies the order $d$ Tur\'an inequalities.  
	 Engel \cite{Engel-2017} showed
	that the overpartition function $\bar{p}(n)$  satisfies the Tur\'an
	inequalities for $n\geq 2$. Liu and Zhang \cite{Liu-Zhang-2021} proved that the
	overpartition function satisfies the higher order Tur\'an inequalities
	for $n\geq 16$. Bringmann, Kane, Rolen and Tripp \cite{Bringmann-Kane-Rolen-Tripp-2021} showed that
	the $k$-colored partition function $p_k(n)$ is log-concave for $n\geq
	6$. Ono, Pujahari and Rolen \cite{Ono-Pujahari-Rolen-2019} proved that the number of the
	MacMahon's plane partitions of $n$ satisfies the order $d\geq 1$ Tur\'an
	inequalities for sufficiently large $n$.

	It should be  noted that the log-concavity of the partition function $p(n)$ for $n\geq 26$ implies that for $n\geq 26$,
	\begin{equation}
		D^2\log p(n-1)<0,
	\end{equation}
	where $D$ is the difference operator with respect to $n$. In general, 
 Chen, Wang and Xie \cite{Chen-Wang-Xie-2016} proved the  following theorem. 
 \begin{theorem}[Chen, Wang and Xie \cite{Chen-Wang-Xie-2016}]
 	For each $r\geq 1$, there exists a positive integer $n(r)$ such that for $n\geq n(r)$, 
 	\begin{equation}\label{eq-diff}
 		(-1)^{r-1} D^{r} \log p(n)>0.
  	\end{equation}
 \end{theorem}
They also noted that  \eqref{eq-diff} is analogous to the positivity of finite differences of the partition function $p(n)$. It says that for given $r\geq 1$,
\begin{equation*}
	D^{r}p(n)>0
\end{equation*}
holds for sufficiently large $n$; see \cite{Gupta-1978}.

	Let us now turn to the broken $k$-diamond partition function and let $\Delta_k(n)$ denote the number of broken $k$-diamond partition of $n$.
	By exploring Sussman's Rademacher-type formula for $\eta$-quotients \cite{Sussman-2017},
	Dong, Ji and Jia \cite{Dong-Ji-Jia-2022} obtained an upper bound and a lower bound for $\Delta_k(n)$ as follows. 
	\begin{theorem}[Dong, Ji and Jia \cite{Dong-Ji-Jia-2022}, Theorem 3.1] \label{Dong-Ji-Jia-the}
	 For $k=1$ or $2$, let $x_k(n)$ be
		\begin{equation}\label{defi-xk}
			x_k(n)=\frac{\pi\sqrt{24n-(2k+2)}}{6},
		\end{equation}
		and 
		\begin{equation}\label{defi-alphak}
	    \alpha_k=\frac{5k+2}{2k+1}.
		\end{equation}
		Then for $k=1$ or $2$  and for $x_k(n) \geq 152 $, or equivalently, for $n\geq 3512$, we have
		\begin{equation*}
		M_k(n)\left(1-\frac{1}{x_k(n)^6}\right)\leq \Delta_k(n)\leq M_k(n)\left(1+\frac{1}{x_k(n)^6}\right).
		\end{equation*}
		where
		\begin{equation}\label{def-Mk}
		M_k(n)=\frac{\alpha_k \pi^3}{18  {x_k(n)^2}}  I_2\left(\sqrt{\alpha_k } {x_k(n)}\right),
		\end{equation}
		and $I_{2}$ is the $2$-nd modified Bessel function of the first kind.
	\end{theorem}
 Based the above theorem, they proved that
	for $k=1$ or $2$, the broken $k$-diamond partition function $\Delta_k(n)$ satisfies the Tur\'an inequalities for $n\geq 1$, which is equivalent to that 
	\begin{equation*}
		D^2\log \Delta_k(n)<0.
	\end{equation*}
	 Furthermore, they derived that for $k=1$ or $2$ and $d\geq 1$, $\{\Delta_k(n)\}_{n\geq 1}$  satisfies the order $d$ Tur\'an inequalities for sufficiently large $n$. As a corollary, they also obtained that 	
	\begin{corollary}
		Let $\Delta_k=\{\Delta_k(n)\}_{n\geq 0}$. For $k=1$ or $2$ and $d\geq 1$, the Jensen polynomial $J_{\Delta_k}^{d,n}(X)$ associated to $\Delta_k$ is hyperbolic for all sufficiently large $n$.
	\end{corollary}
	
 By Theorem \ref{Dong-Ji-Jia-the}, one can easily obtain the following corollary.
	\begin{corollary}\label{cor-Dong-Ji-Jia}
	Let $x_k(n)$ be defined as \eqref{defi-xk}.  
	For $k=1$ or $2$, let 	
	\begin{equation}\label{def-lambda}
		\Lambda_k(n)=\frac{M_k(n-1)M_k(n+1)}{M_k(n)^2},
	\end{equation}
	and
	\begin{equation}\label{def-Theta}
		\Theta_k(n)= \frac{\Delta_k(n-1)\Delta_k(n+1)}{\Delta_k(n)^2},
	\end{equation}
	where $M_k(n)$ is defined as in \eqref{def-Mk}.
	Then for  $k=1$ or $2$ and for $x_k(n)\geq 152$, or equivalently, $n\geq 3512$,
	\begin{equation}\label{eq-delta-1}
		\Lambda_k(n) g_k(n) \leq \Theta_k(n) \leq \Lambda_k(n) G_k(n),
	\end{equation}
	where 
	\begin{equation}\label{def-g-k-n}
		g_k(n)=\frac{\left(1-\frac{1}{x_k(n-1)^6}\right)\left(1-\frac{1}{x_k(n+1)^6}\right)}{\left(1+\frac{1}{x_k(n)^6}\right)^2},
	\end{equation}
	and 
	\begin{equation}\label{def-G-k-n}
		G_k(n)=\frac{\left(1+\frac{1}{x_k(n-1)^6}\right)\left(1+\frac{1}{x_k(n+1)^6}\right)}{\left(1-\frac{1}{x_k(n)^6}\right)^2}.
	\end{equation}
	\end{corollary}
	
	In this paper, we will further study this corollary and prove that the broken $k$-diamond partition function $\Delta_k(n)$ satisfies the higher order Tur\'an inequalities when $n\geq 6$, where $k=1$ or $2$.
	To this end, we first derive an efficient lower and upper bounds for $\Theta_k(n)$ as follows.
	
		\begin{theorem}\label{the-Theta-bounding}
	Let $\Theta_k(n)$ be defined as in \eqref{def-Theta}.
		For  $k=1$ or $2$ and for $x_k(n)\geq 315$, or equivalently, for $n\geq 15081$,  we have 
		\begin{align}
		\Theta_k(n)> 
		 1-\frac{\pi ^4 \sqrt{\alpha_k}}{9 x_k(n)^3}+\frac{5\pi ^4}{9 x_k(n)^4}-\frac{5 \pi ^4}{8 x_k(n)^5 \sqrt{\alpha_k}}+\frac{-\frac{300}{\alpha_k^3}-10-\frac{5 \pi ^4}{6 \alpha_k}}{x_k(n)^6},
		\end{align}
		and 
		\begin{align}\label{eq-theta-bounding-upper}
			\Theta_k(n)<
		1-\frac{\pi ^4 \sqrt{\alpha_k}}{9 x_k(n)^3}+\frac{5\pi ^4}{9 x_k(n)^4}-\frac{5 \pi ^4}{8 x_k(n)^5 \sqrt{\alpha_k}}+\frac{\frac{ \pi ^8 \alpha_k}{81}+\frac{292}{\alpha_k^3}+5}{x_k(n)^6}.
		\end{align}
	\end{theorem}
	
	It should be noted that exploring the upper bound for $\Theta_k(n)$ in \eqref{eq-theta-bounding-upper}, one can also prove that for $k=1$ or $2$, the broken $k$-diamond partition function is log-concave.
	
	Employing these bounds for $\Theta_k(n)$ in Theorem \ref{the-Theta-bounding}, we  show that 
	
	\begin{theorem}\label{the-Theta-mono}
For $n \geq 5$,
	\begin{equation}
		\Theta_1(n)< \Theta_1(n+1),
	\end{equation}
	and 
	for $n\geq 7$,
	\begin{equation}
		\Theta_2(n)<\Theta_2(n+1).
	\end{equation}
%
	\end{theorem}

The above theorem  can be  restated as follows.
	\begin{theorem}\label{the-Theta-D}
	For $n \geq 5$,
	\begin{equation}
		D^3\log\Delta_1(n-1)>0,
	\end{equation}
	and 
	for $n\geq 7$,
	\begin{equation}
		D^3\log\Delta_2(n-1)>0.
	\end{equation}
		where $D$ is the difference operator with respect to $n$.
	\end{theorem}
	
In general, we propose the following conjecture, which is analogous to the positivity of finite differences of the logarithm of the partition function $p(n)$.
\begin{conjecture}
Let $\Delta_k(n)$ be the number of the  $k$-diamond broken partition of $n$. For any $k\geq 1$ and $r\geq 1$, there exists a positive integer $n_k(r)$ such that for $n\geq n_k(r)$,
	\begin{equation}
		(-1)^{r-1}D^r \log \Delta_k(n)>0.
	\end{equation}
\end{conjecture}
	
	
   With Theorem \ref{the-Theta-bounding} and Theorem \ref{the-Theta-mono} in hand, we can prove that for $k=1$ or $2$, 
   $\{\Delta_k(n)\}_{n\geq 6}$ satisfies the higher order  Tur\'an inequalities, which implies that for $k=1$ or $2$, the Jensen polynomial $J_{\Delta_k}^{3,n}(X)$ has only real zeros when $n\geq 6$.
    
   	\begin{theorem}
		For $k=1$ or $2$ and $n\geq 6$, we have 
		\begin{align}
			\nonumber 
			&4\big(\Delta_k(n)^2-\Delta_k(n-1)\Delta_k(n+1)\big)\big(\Delta_k(n+1)^2-\Delta_k(n)\Delta_k(n+2)\big)\\
			&\qquad -\big(\Delta_k(n)\Delta_{k}(n+1)-\Delta_k(n-1)\Delta_k(n+2)\big)^2>0.
		\end{align}
	\end{theorem} 
In fact, we  shall prove the following equivalent form.
\begin{theorem}\label{the-Theta-k}
	Let $\Theta_k(n)$ be defined as in \eqref{def-Theta}.
	For $k=1$ or $2$ and for $n\geq 6$, we have
	\begin{equation}
		4(1-\Theta_k(n))(1-\Theta_k(n+1))-(1-\Theta_k(n)\Theta_k(n+1))>0.
	\end{equation}		
	\end{theorem}

\section{AN INEQUALITY ON THE BESSEL FUNCTION}

In this section, we shall  bound the Bessel function of the first kind effectively, which will play a key role in our paper. In fact, we prove the following inequalities.
 	
 	\begin{theorem}
 	Let $\nu \geq 2$ and let $I_{\nu}(z)$ be the $\nu$-th modified Bessel function of the first kind. 
 	Then we have that for $\nu \geq 2$ and $z\geq \frac{\left(\nu+\frac{11}{2}\right)^6}{120}$,
 		\begin{align}
    		\nonumber &\left|I_{\nu}(z)e^{-z}\sqrt{2 \pi z}-\sum _{i=0}^5 \frac{\prod _{j=1}^i \left(-j+\nu -\frac{1}{2}+1\right) \left(j+\nu +\frac{1}{2}-1\right)}{i! (-2 z)^i}\right|\\[9pt]
    		\nonumber & \quad \leq 
     		 \frac{52e^{-z} }{17\Gamma(\nu+\frac{1}{2})} 
			\sum_{i=0}^{5}\left| \binom{\nu-\frac{1}{2}}{i} \right| \frac{z^{\nu-\frac{1}{2}}}{2^i}
			+\frac{e^{-z}}{2^{\nu-\frac{1}{2}}\Gamma(\nu+\frac{1}{2})} z^{\nu+\frac{1}{2}}\\[9pt]
			 &\qquad \quad   +\left| \frac{ \left(\nu^2 -\frac{121}{4}\right) \left(\nu^2 -\frac{81}{4}\right) \left(\nu^2 -\frac{49}{4}\right) \left(\nu^2 -\frac{25}{4}\right) \left(\nu^2 -\frac{9}{4}\right) \left(\nu^2 -\frac{1}{4}\right)}{6! 2^{\nu-\frac{1}{2}}z^6}\right|\max\left\{2^{\nu-\frac{13}{2}},1\right\}. \label{eq-I-nu}
		\end{align}
	In particular, for  $\nu=2$,  we have that for $z\geq  \frac{\left(\frac{15}{2}\right)^6}{120}\approx 1483.2$,
	\begin{align}\label{eq-I-2}
	 \left|I_{2}(z)e^{-z}\sqrt{2 \pi z}-1+\frac{15}{8 z}-\frac{105}{128 z^2}-\frac{315}{1024 z^3}-\frac{10395}{32768 z^4}-\frac{135135}{262144 z^5}\right|\leq \frac{73}{z^6}. 
	 	\end{align}
 	\end{theorem}
 
 	{\it Proof:} We begin with the following integral representation \cite{Waston-1944} of the $\nu$-th modified Bessel function:
	\begin{equation}
		I_{\nu}(z)=\frac{\left(\frac{z}{2}\right)^{\nu}}{\Gamma(\nu+\frac{1}{2})\sqrt{\pi}} \int_{-1}^{1} (1-t^2)^{\nu-\frac{1}{2}} e^{z t}{\rm d t}.
	\end{equation}
	Then we have that  
		\begin{align}
			\nonumber I_{\nu}(z)&=\frac{\left(\frac{z}{2}\right)^{\nu}}{\Gamma(\nu+\frac{1}{2})} \int_{-1}^{1} (1-t^2)^{\nu-\frac{1}{2}} e^{z t}{\rm d t}\\[9pt]
			&=\frac{\left(\frac{z}{2}\right)^{\nu}}{\Gamma(\nu+\frac{1}{2})} 
				\int_{0}^{1} (1-t^2)^{\nu-\frac{1}{2}} e^{z t}{\rm d t}
				+\frac{\left(\frac{z}{2}\right)^{\nu}}{\Gamma(\nu+\frac{1}{2})} 
				\int_{-1}^{0} (1-t^2)^{\nu-\frac{1}{2}} e^{z t}{\rm d t}.
		\end{align}
		
	Setting $1-t=u$ in the first integral, we get that 
	\begin{equation}\label{eq-I-nu-2}
		I_{\nu}(z)=\frac{\left(\frac{z}{2}\right)^{\nu}}{\Gamma(\nu+\frac{1}{2})} 
				\int_{0}^{1} (2-u)^{\nu-\frac{1}{2}}u^{\nu-\frac{1}{2}} e^{z(1-u)}{\rm du}
					+\frac{\left(\frac{z}{2}\right)^{\nu}}{\Gamma(\nu+\frac{1}{2})} 
				\int_{-1}^{0} (1-t^2)^{\nu-\frac{1}{2}} e^{z t}{\rm d t}.
	\end{equation}
		
	By Taylor's Theorem, we have  
		\begin{align}
			\nonumber (2-u)^{\nu-\frac{1}{2}}&= 2^{\nu -\frac{1}{2}}
				- \left(\nu -\frac{1}{2}\right) 2^{\nu -\frac{3}{2}} u
				+\frac{1}{2!}\left(\nu -\frac{1}{2}\right) \left(\nu -\frac{3}{2}\right)2^{\nu -\frac{5}{2}}u^2\\
			\nonumber &\qquad 
				-\frac{1}{3!}  \left(\nu -\frac{1}{2}\right)  \left(\nu -\frac{3}{2}\right) \left(\nu -\frac{5}{2}\right)2^{\nu -\frac{7}{2}}  u^3\\
			\nonumber &\qquad 
				+\frac{1}{4!} \left(\nu -\frac{1}{2}\right) \left(\nu -\frac{3}{2}\right) \left(\nu -\frac{5}{2}\right) \left(\nu -\frac{7}{2}\right) 2^{\nu -\frac{9}{2}} u^4\\
			\nonumber &\qquad	-\frac{1}{5!}    \left(\nu -\frac{1}{2}\right) \left(\nu -\frac{3}{2}\right)  \left(\nu -\frac{5}{2}\right) \left(\nu -\frac{7}{2}\right) \left(\nu -\frac{9}{2}\right)2^{\nu -\frac{11}{2}} u^5\\
			\nonumber &\qquad +C_{\nu}(u)u^6, \\
			&=\sum_{i=0}^{5} \binom{\nu-\frac{1}{2}}{i} (-1)^i 2^{\nu-\frac{1}{2}-i} u^{i}+C_{\nu}(u)u^6,
			\label{eq-bes-int-1}
		\end{align}
	where for some $\xi \in [0,1]$,
	\begin{align*}
		C_{\nu}(u)&=\frac{1}{6!} \left[\frac{\partial^6}{\partial u^6}(2-u)^{\nu-\frac{1}{2}}\right]_{u=\xi}\\
		&=\frac{1}{6!}   \left(\nu -\frac{1}{2}\right) \left(\nu -\frac{3}{2}\right)
		\left(\nu -\frac{5}{2}\right) 
		 \left(\nu -\frac{7}{2}\right)\left(\nu -\frac{9}{2}\right)  \left(\nu -\frac{11}{2}\right)(2-\xi)^{\nu-\frac{13}{2}}.
	\end{align*}
	And we note that for $\nu\geq 2$,
	\begin{align}\label{ineq-C-nu}
		|C_{\nu}(u)|\leq \frac{1}{6!}   \left(\nu -\frac{1}{2}\right) \left(\nu -\frac{3}{2}\right)
		\left|\nu -\frac{5}{2}\right|
		 \left|\nu -\frac{7}{2}\right|\left|\nu -\frac{9}{2}\right|  \left|\nu -\frac{11}{2}\right| \max\left\{2^{\nu-\frac{13}{2}},1\right\}.
	\end{align}

 Plugging \eqref{eq-bes-int-1} in \eqref{eq-I-nu-2}, we get that 
	\begin{align}
		\nonumber	I_{\nu}(z)&=\frac{\left(\frac{z}{2}\right)^{\nu}}		{				\Gamma(\nu+\frac{1}{2})\sqrt{\pi}} 
						\int_{0}^{1} (2-u)^{\nu-\frac{1}{2}}u^{\nu-\frac{1}{2}} e^{z(1-u)}{\rm du}
					  +\frac{\left(\frac{z}{2}\right)^{\nu}}{\Gamma(\nu+\frac{1}{2})\sqrt{\pi}} 
					\int_{-1}^{0} (1-t^2)^{\nu-\frac{1}{2}} e^{z t}{\rm d t}\\[9pt]
		\nonumber &=\frac{\left(\frac{z}{2}\right)^{\nu}}{\Gamma(\nu+\frac{1}			{2})\sqrt{\pi}} 
						\int_{0}^{1}\left[\sum_{i=0}^{5} \binom{\nu-\frac{1}{2}}{i} (-1)^i 2^{\nu-\frac{1}{2}-i} u^{i}+C_{\nu}(u)u^6\right]
				u^{\nu-\frac{1}{2}} e^{z(1-u)}{\rm du}\\[9pt]
		\nonumber &\qquad +\frac{\left(\frac{z}{2}\right)^{\nu}}{\Gamma(\nu+				\frac{1}{2})\sqrt{\pi}} 
						\int_{-1}^{0} (1-t^2)^{\nu-\frac{1}{2}} e^{z t}{\rm d t}\\[9pt]
		\nonumber &=\frac{\left(\frac{z}{2}\right)^{\nu}}{\Gamma(\nu+\frac{1}			{2})\sqrt{\pi}} 
					\int_{0}^{1} \sum_{i=0}^{5} \binom{\nu-\frac{1}{2}}{i} (-1)^i 2^{\nu-\frac{1}{2}-i}
						u^{\nu-\frac{1}{2}+i} e^{z(1-u)}{\rm du}\\[9pt]
	   \nonumber  &\qquad 	+\frac{\left(\frac{z}{2}\right)^{\nu}}{\Gamma(\nu+				\frac{1}{2})\sqrt{\pi}} 
					\int_{0}^{1}C_{\nu}(u)u^{\nu+\frac{11}{2}} e^{z(1-u)}{\rm du}+\frac{\left(\frac{z}{2}\right)^{\nu}}{\Gamma(\nu+\frac{1}{2})\sqrt{\pi}} 
						\int_{-1}^{0} (1-t^2)^{\nu-\frac{1}{2}} e^{z t}{\rm d t}\\[9pt]
		\nonumber &=\frac{\left(\frac{z}{2}\right)^{\nu}}{\Gamma(\nu+\frac{1}				{2})\sqrt{\pi}} 
						\left(\int_{0}^{\infty}-\int_{1}^{\infty}\right)\sum_{i=0}^{5} \binom{\nu-\frac{1}{2}}{i} (-1)^i 2^{\nu-\frac{1}{2}-i}
						u^{\nu-\frac{1}{2}+i} e^{z(1-u)}{\rm du}\\[9pt]
	   \nonumber  &\qquad 	+\frac{\left(\frac{z}{2}\right)^{\nu}}{\Gamma(\nu+				\frac{1}{2})\sqrt{\pi}} 
					\int_{0}^{1}C_{\nu}(u)u^{\nu+\frac{11}{2}} e^{z(1-u)}{\rm du}+\frac{\left(\frac{z}{2}\right)^{\nu}}{\Gamma(\nu+\frac{1}{2})\sqrt{\pi}} 
						\int_{-1}^{0} (1-t^2)^{\nu-\frac{1}{2}} e^{z t}{\rm d t}\\[9pt]
		\nonumber &=\frac{e^z}{\Gamma(\nu+\frac{1}{2})\sqrt{2 \pi z}} 
					\sum_{i=0}^{5} \binom{\nu-\frac{1}{2}}{i} \frac{1}{(-2z)^i}
					\left(\int_{0}^{\infty}-\int_{1}^{\infty}\right)
					(zu)^{\nu-\frac{1}{2}+i} e^{-zu}\cdot z{\rm du}\\[9pt]
	   \nonumber  &\qquad 	+\frac{e^z}{\Gamma(\nu+\frac{1}{2} )\sqrt{2 \pi z}}
	    				\cdot 
	   				 \frac{C_{\nu}(u)}{2^{\nu-\frac{1}{2}}z^6}
					\left(\int_{0}^{\infty}-\int_{1}^{\infty}\right)(zu)^{\nu+\frac{11}{2}} e^{-zu}\cdot z{\rm du}\\[9pt]
		\nonumber &\qquad +\frac{\left(\frac{z}{2}\right)^{\nu}}{\Gamma(\nu+			\frac{1}{2})\sqrt{\pi}} 
					\int_{-1}^{0} (1-t^2)^{\nu-\frac{1}{2}} e^{z t}{\rm d t}\\[9pt]
		\nonumber &=\frac{e^z}{\Gamma(\nu+\frac{1}{2})\sqrt{2 \pi z}} 
					\sum_{i=0}^{5} \binom{\nu-\frac{1}{2}}{i} \frac{1}{(-2z)^i}
						\int_{0}^{\infty} t^{\nu+i+\frac{1}{2}-1} e^{-t}{\rm dt}\\[9pt]
		\nonumber	&\qquad -\frac{e^z}{\Gamma(\nu+\frac{1}{2})\sqrt{2 \pi z}} 
						\sum_{i=0}^{5} \binom{\nu-\frac{1}{2}}{i} \frac{1}{(-2z)^i}\int_{z}^{\infty}
						t^{\nu+i+\frac{1}{2}-1} e^{-t}{\rm dt}\\[9pt]
	    \nonumber	&\qquad 	 	+\frac{e^z}{\Gamma(\nu+\frac{1}{2} )\sqrt{2 			\pi z}}
	  					  \cdot 
	  					  \frac{C_{\nu}(u)}{2^{\nu-\frac{1}{2}}z^6}
							\left(\int_{0}^{\infty}-\int_{z}^{\infty}\right)t^{\nu+\frac{11}{2}} e^{-t}{\rm dt}\\[9pt]
		\nonumber	&\qquad +\frac{\left(\frac{z}{2}\right)^{\nu}}{\Gamma(\nu+			\frac{1}{2})\sqrt{\pi}} 
					\int_{-1}^{0} (1-t^2)^{\nu-\frac{1}{2}} e^{z t}{\rm d t}\\[9pt]
		\nonumber	&=\frac{e^z}{\Gamma(\nu+\frac{1}{2})\sqrt{2 \pi z}} 
					\sum_{i=0}^{5} \binom{\nu-\frac{1}{2}}{i} \frac{1}{(-2z)^i} \Gamma\left(\nu+i+\frac{1}{2}\right)\\[9pt]
		\nonumber	&\qquad -\frac{e^z}{\Gamma(\nu+\frac{1}{2})\sqrt{2 \pi z}} 
						\sum_{i=0}^{5} \binom{\nu-\frac{1}{2}}{i} \frac{1}{(-2z)^i}
						\Gamma\left(\nu+i+\frac{1}{2},z\right)\\[9pt]
		\nonumber	&\qquad 	+\frac{e^z}{\Gamma(\nu+\frac{1}{2} )\sqrt{2 \pi 			z}}
	   					 \cdot 
	  					  \frac{C_{\nu}(u)}{2^{\nu-\frac{1}{2}}z^6}
	  						  \Bigg(\Gamma\left(\nu+\frac{13}{2}\right)-\Gamma\left(\nu+\frac{13}{2},z\right)\Bigg)\\[9pt]
		\nonumber	&\qquad +\frac{\left(\frac{z}{2}\right)^{\nu}}{\Gamma(\nu+				\frac{1}{2})\sqrt{\pi}} 
						\int_{-1}^{0} (1-t^2)^{\nu-\frac{1}{2}} e^{z t}{\rm d t}\\[9pt]
		\nonumber	&=\frac{e^z}{\sqrt{2 \pi z}} 
						\sum _{i=0}^5 \frac{\prod _{j=1}^i \left(-j+\nu -\frac{1}{2}+1\right) \left(j+\nu +\frac{1}{2}-1\right)}{i! (-2 z)^i}\\[9pt]
		\nonumber	&\qquad -\frac{e^z}{\Gamma(\nu+\frac{1}{2})\sqrt{2 \pi z}} 
					\sum_{i=0}^{5} \binom{\nu-\frac{1}{2}}{i} \frac{1}{(-2z)^i}
					\Gamma\left(\nu+i+\frac{1}{2},z\right)\\[9pt]
		\nonumber	&\qquad 	+\frac{e^z}{\Gamma(\nu+\frac{1}{2} )\sqrt{2 \pi 			z}}
	   					 \cdot 
	   					 \frac{C_{\nu}(u)}{2^{\nu-\frac{1}{2}}z^6}
	  					  \Bigg(\Gamma\left(\nu+\frac{13}{2}\right)-\Gamma\left(\nu+\frac{13}{2},z\right)\Bigg)\\[9pt]
		&\qquad +\frac{\left(\frac{z}{2}\right)^{\nu}}{\Gamma(\nu+\frac{1}{2})\sqrt{\pi}} 
			\int_{-1}^{0} (1-t^2)^{\nu-\frac{1}{2}} e^{z t}{\rm d t}. \label{eq-I-nu-main}
	\end{align}
	Thus we have that 
	\begin{align}
		\nonumber &I_{\nu}(z)e^{-z}\sqrt{2 \pi z}-\sum _{i=0}^5 \frac{\prod _{j=1}^i \left(-j+\nu -\frac{1}{2}+1\right) \left(j+\nu +\frac{1}{2}-1\right)}{i! (-2 z)^i}\\
		\nonumber	&\quad = -\frac{1}{\Gamma(\nu+\frac{1}{2})} 
					\sum_{i=0}^{5} \binom{\nu-\frac{1}{2}}{i} \frac{1}{(-2z)^i}
					\Gamma\left(\nu+i+\frac{1}{2},z\right)\\[9pt]
		\nonumber	&\qquad  \quad	+\frac{1}{\Gamma(\nu+\frac{1}{2})}
	   					 \cdot 
	   					 \frac{C_{\nu}(u)}{2^{\nu-\frac{1}{2}}z^6}
	  					  \Bigg(\Gamma\left(\nu+\frac{13}{2}\right)-\Gamma\left(\nu+\frac{13}{2},z\right)\Bigg)\\[9pt]
		&\qquad \quad +\frac{z^{\nu+\frac{1}{2}}e^{-z}}{2^{\nu-\frac{1}{2}}\Gamma(\nu+\frac{1}{2})} 
			\int_{-1}^{0} (1-t^2)^{\nu-\frac{1}{2}} e^{z t}{\rm d t}. \label{eq-I-nu-4}
	\end{align}
	Thus to prove \eqref{eq-I-nu}, we shall estimate the three terms in \eqref{eq-I-nu-4}.
	Clearly, the  integral term in \eqref{eq-I-nu-4} can be easily bounded as follows: 
	\begin{align}\label{ineq-I-nu-3}
      	\left|\frac{z^{\nu+\frac{1}{2}}e^{-z}}{2^{\nu-\frac{1}{2}}\Gamma(\nu+\frac{1}{2})} 
			\int_{-1}^{0} (1-t^2)^{\nu-\frac{1}{2}} e^{z t}{\rm d t}\right|
		\leq \frac{e^{-z}}{2^{\nu-\frac{1}{2}}\Gamma(\nu+\frac{1}{2})}
		z^{\nu+\frac{1}{2}},
   \end{align}
   and by applying \eqref{ineq-C-nu}, we find that for $\nu \geq 2$, the second term in \eqref{eq-I-nu-4} can be bounded as follows:  
		\begin{align}
        	\nonumber &\left|	\frac{1}{\Gamma(\nu+\frac{1}{2} )}
	    		 \cdot 
	     		\frac{C_{\nu}(u)}{2^{\nu-\frac{1}{2}}z^6}
	    		 \Bigg(\Gamma\left(\nu+\frac{13}{2}\right)-\Gamma\left(\nu+\frac{13}{2},z\right)\Bigg)\right|\\[9pt]
	    	\nonumber &\qquad \leq \frac{1}{\Gamma(\nu+\frac{1}{2} )}
	   			 \cdot 
	   			  \frac{1}{2^{\nu-\frac{1}{2}}z^6}\cdot
	    		 \Gamma\left(\nu+\frac{13}{2}\right)\cdot\left|C_{\nu}(u)\right|
	    			 \\[9pt]
	   		\nonumber &\qquad =\frac{1}{2^{\nu-\frac{1}{2}}z^6}\cdot
	     		\left|C_{\nu}(z) \left(\nu +\frac{11}{2}\right) \left(\nu +\frac{9}{2}\right) \left(\nu +\frac{7}{2}\right) \left(\nu +\frac{5}{2}\right) \left(\nu +\frac{3}{2}\right) \left(\nu +\frac{1}{2}\right)\right|
	    		 \\[9pt]
	  		&\qquad \leq
	    	 	\left| \frac{ \left(\nu^2 -\frac{121}{4}\right) \left(\nu^2 -\frac{81}{4}\right) \left(\nu^2 -\frac{49}{4}\right) \left(\nu^2 -\frac{25}{4}\right) \left(\nu^2 -\frac{9}{4}\right) \left(\nu^2 -\frac{1}{4}\right)
						}{6! 2^{\nu-\frac{1}{2}}z^6}\right|\max\left\{2^{\nu-\frac{13}{2}},1\right\}. \label{ineq-I-nu-2}
		\end{align}
	To bound the first term in \eqref{eq-I-nu-4}, we note that for $z\geq \frac{a^6}{120}$ and $a\geq \frac{5}{2}$, 
	\begin{equation}\label{ineq-inco-gamma}
		\Gamma(a,x)\leq \frac{52}{17}z^{a-1}e^{-z};
	\end{equation}
	see \cite{Bringmann-Kane-Rolen-Tripp-2021}.
	Then by \eqref{ineq-inco-gamma}, we obtain that for $\nu \geq 2$ and $z\geq \frac{\left(\nu+\frac{11}{2}\right)^6}{120}$,
	\begin{align}
		\nonumber &\left|- \frac{1}{\Gamma(\nu+\frac{1}{2})} 
					\sum_{i=0}^{5} \binom{\nu-\frac{1}{2}}{i} \frac{1}{(-2z)^i}
					\Gamma\left(\nu+i+\frac{1}{2},z\right)\right|\\[9pt]
		\nonumber &\qquad \leq \frac{1}{\Gamma(\nu+\frac{1}{2})} 
					\sum_{i=0}^{5}\left| \binom{\nu-\frac{1}{2}}{i} \right| \frac{1}{(2z)^i}\cdot \frac{52}{17} z^{\nu+i+\frac{1}{2}-1} e^{-z}\\[9pt]
		&\qquad = \frac{52 e^{-z}}{17\Gamma(\nu+\frac{1}{2})} 
		\sum_{i=0}^{5}\left| \binom{\nu-\frac{1}{2}}{i} \right| \frac{z^{\nu-\frac{1}{2}}}{2^i}. \label{ineq-I-nu-1}
	\end{align}
	Combing \eqref{eq-I-nu-4}, \eqref{ineq-I-nu-3}, \eqref{ineq-I-nu-2} and \eqref{ineq-I-nu-1}, we are led to \eqref{eq-I-nu}.
	When $\nu=2$,  we have that for $z\geq \frac{\left(\frac{15}{2}\right)^6}{120}  \approx 1483.2$,
		\begin{align*}
			 &\left|I_{2}(z)e^{-z}\sqrt{2 \pi z}-1+\frac{15}{8 z}-\frac{105}{128 z^2}-\frac{315}{1024 z^3}-\frac{10395}{32768 z^4}-\frac{135135}{262144 z^5}\right|\\[9pt]
    		& \quad \quad \leq 
    		  \frac{65793 e^{-z}}{8704 \sqrt{\pi }}  z^{3/2}
				+\frac{2e^{-z}}{3\sqrt{2\pi}} z^{\frac{5}{2}}
				+\frac{4729725}{65536}z^{-6}\\
			&\quad \quad \leq 5 e^{-z}z^{5/2}+\frac{4729725}{65536}z^{-6}.
		\end{align*}
	
	We now claim that for $z\geq 31$,
		\begin{align*}
			e^{-z}z^{5/2}\leq \frac{54403}{327680} z^{-6}.
		\end{align*}
	Let 
		\begin{align*}
		F(t)=e^{-z}z^{17/2}.
		\end{align*}
	Then
		\[F'(t)=\frac{-2 t+17}{2} e^{-t} t^{15/2}.\]
	For $t> \frac{17}{2}$, we see $F'(t)$ is negative. Then $F(t)$ is decreasing on $(\frac{17}{2},\infty)$. Thus for $z\geq 31$, we have
		\begin{align*}
			e^{-z}z^{17/2}\leq F(31)\approx 0.163 <\frac{54403}{327680} \approx 0.166,
		\end{align*}
	which confirms our claim. Therefore, we obtain that for $z\geq \frac{\left(\frac{15}{2}\right)^6}{120}  \approx 1483.2$,
		\begin{align*}
	 		&\left|I_{2}(z)e^{-z}\sqrt{2 \pi z}-1+\frac{15}{8 z}-\frac{105}{128 z^2}-\frac{315}{1024 z^3}-\frac{10395}{32768 z^4}-\frac{135135}{262144 z^5}\right|\\[9pt]
			&\qquad \leq 5 e^{-z}z^{5/2}+\frac{4729725}{65536}z^{-6}\\[9pt]
			&\qquad \leq \frac{54403}{65536} z^{-6}+\frac{4729725}{65536}z^{-6}=\frac{73}{z^6}.
		\end{align*}

	This completes the proof. 
	\qed

\section{BOUNDING $\Lambda_k(n)$}\label{sec-Lambda}

In this section, we shall give an upper bound and a lower bound for 
\[\Lambda_k(n)=\frac{M_k(n-1)M_{k}(n+1)}{M_k(n)^2},\]
as defined in \eqref{def-lambda},
where $M_k(n)$ is defined as in \eqref{def-Mk}: 
\[M_k(n)=\frac{\alpha_k \pi^3}{18  {x^2_k(n)}}  I_2\left(\sqrt{\alpha_k } {x_k(n)}\right),\]
and $x_k(n)$ is defined as in \eqref{defi-xk}:
\begin{equation*}
	x_k(n)=\frac{\pi\sqrt{24n-(2k+2)}}{6},
\end{equation*}
and $\alpha_k=\frac{5k+2}{2k+1}$. 
In fact, we prove the following theorem.

\begin{theorem}\label{the-Lambda-bounding}
Let $\Lambda_k(n)$ be  defined as  in \eqref{def-lambda}.
For $k=1$ or $2$ and for $n \geq 2$, we have 
	\begin{align}
	\nonumber \Lambda_k(n) &< 
	\left(1+\frac{5\pi ^4}{9 x_k(n)^4}+\frac{ \pi ^8}{3 x_k(n)^8}\right)\left(1-\frac{\sqrt{\alpha_k}\pi^4}{9x_k(n)^3}+\frac{\alpha_k\pi^8}{81x_k(n)^6}\right)\\
	\label{eq-Lambda-leq} &\qquad \times \left( 1-\frac{5 \pi ^4}{8 \sqrt{\alpha_k} x_k(n)^5}+\frac{292}{\alpha_k^3 x_k(n)^6}\right),
\end{align}
and 
\begin{align}
	\nonumber \Lambda_k(n)&> \left(1+\frac{5\pi ^4}{9 x_k(n)^4}+\frac{ 5\pi ^8}{18 x_k(n)^8}\right)\left(1-\frac{\sqrt{\alpha_k}\pi^4}{9x_k(n)^3}-\frac{5\sqrt{\alpha_k}\pi^8}{162x_k(n)^7}\right)\\
	\label{eq-Lambda-geq}	&\qquad \times\left( 1-\frac{5 \pi ^4}{8 \sqrt{\alpha_k} x_k(n)^5}-\frac{5 \pi ^4}{6 \alpha_k x_k(n)^6}-\frac{300}{\alpha_k^3 x_k(n)^6}\right).
\end{align}
\end{theorem}

The proof of this theorem is based on Lemma \ref{lem-H-bounding}, which gives tow  inequalities that involves the second modified Bessel function of the first kind. 

\begin{lemma}\label{lem-H-bounding}
 Let
		\begin{align}\label{eq-H-k}
			H_k(n)=	\frac{I_2\left(\sqrt{\alpha_k} x_k(n-1)\right) I_2\left(\sqrt{\alpha_k} x_k(n+1)\right)}{I_2\left(\sqrt{\alpha_k} x_k(n)\right)^2}.
		\end{align}
 Then we have that for $x_k(n-1)\geq 971$, or equivalently, for $n\geq 143296$,
 \begin{align}
	\nonumber	H_k(n)&<\frac{x_k(n)}{\sqrt{x_k(n-1) x_k(n+1)}}
						\left(1- \frac{\sqrt{\alpha_k}\pi^4}{9x_k(n)^3}+\frac{\alpha_k\pi^8}{81x_k(n)^6}\right)\\[6pt]
    \label{eq-H-Leq}  &\qquad \times 
    					\left( 1-\frac{5 \pi ^4}{8 			  \sqrt{\alpha_k} x_k(n)^5}+\frac{292}{\alpha_k^3 x_k(n)^6}\right),\\[9pt]
	\nonumber	H_k(n)&> \frac{x_k(n)}{\sqrt{x_k(n-1) x_k(n+1)}}
						\left(1-\frac{\sqrt{\alpha_k}\pi^4}{9x_k(n)^3}-\frac{5\sqrt{\alpha_k}\pi^8}{162x_k(n)^7}\right)
						\\[6pt]
	\label{eq-H-Geq}  &\qquad \times\left( 1-\frac{5 \pi ^4}{8 \sqrt{\alpha_k} x_k(n)^5}-\frac{5 \pi ^4}{6 \alpha_k x_k(n)^6}-\frac{300}{\alpha_k^3 x_k(n)^6}\right).
 \end{align}
\end{lemma}

{\it Proof:}
For convenience,  set
\begin{equation*}
\begin{aligned}
	\gamma_1=\frac{15}{8 \sqrt{\alpha_k}}, 
	\quad 
	\gamma_2=\frac{105}{128 \alpha_k},
	\quad 
	\gamma _3=\frac{315}{1024 \alpha_k^{3/2}},\\[9pt]
    \end{aligned}
\end{equation*}
\begin{equation*}
\begin{aligned}
    \gamma_4=\frac{10395}{32768 \alpha_k^2},
    \quad
    \gamma_5=\frac{135135}{262144 \alpha_k^{5/2}},
    \quad 
    \gamma_6=\frac{73}{\alpha_k^3},
\end{aligned}
\end{equation*}

and let
\begin{align*}
	\varphi_k(t)=1-\frac{\gamma _1}{t}+\frac{\gamma _2}{t^2}+\frac{\gamma _3}{t^3}+\frac{\gamma _4}{t^4}+\frac{\gamma _5}{t^5}+\frac{\gamma _6}{t^6},\\[9pt]
	\varPhi_k(t)=1-\frac{\gamma _1}{t}+\frac{\gamma _2}{t^2}+\frac{\gamma _3}{t^3}+\frac{\gamma _4}{t^4}+\frac{\gamma _5}{t^5}-\frac{\gamma _6}{t^6}.
\end{align*}
For $k=1$ or $2$, by \eqref{eq-I-2}, we find that for $t \geq \max\{\frac{\left(\frac{15}{2}\right)^6}{120\sqrt{\alpha_1}}, 
\frac{\left(\frac{15}{2}\right)^6}{120\sqrt{\alpha_2}}\}\approx 970.95$,
\begin{equation}\label{eq-I-alpha-t}
	\frac{e^{\sqrt{\alpha_k} t}}{\sqrt{2\pi \sqrt{\alpha_k}t}}\varphi_k(t)\leq I_2(\sqrt{\alpha_k}t)\leq \frac{e^{\sqrt{\alpha_k} t}}{\sqrt{2\pi \sqrt{\alpha_k}t}}\varPhi_k(t).
\end{equation}
Applying \eqref{eq-I-alpha-t} to $H_k(n)$ (see \eqref{eq-H-k}), we obtain that for $x_k(n-1)\geq 971$, or equivalently, for $n\geq 143296$,
\begin{align}
	\nonumber 
	H_k(n)&\leq \frac{\frac{e^{\sqrt{\alpha_k}x_k(n-1)}}{\sqrt{2\pi \sqrt{\alpha_k} x_k(n-1)}} \varPhi_k(x_k(n-1)) \cdot \frac{e^{\sqrt{\alpha_k}x_k(n+1)}}{\sqrt{2\pi \sqrt{\alpha_k} x_k(n+1)}} \varPhi_k(x_k(n+1))}
	{\frac{e^{2\sqrt{\alpha_k}x_k(n)}}{2\pi \sqrt{\alpha_k} x_k(n)} \varphi_k(x_k(n))^2}\\[9pt]
	&= \frac{x_k(n)}{\sqrt{x_k(n-1) x_k(n+1)}}
	e^{\sqrt{\alpha_k } (x_k(n-1)+x_k(n+1)-2x_k(n))}L_k(n), \label{eq-H-leq-L}
\end{align}
and 
\begin{align}
	\nonumber 
	H_k(n)&\geq \frac{\frac{e^{\sqrt{\alpha_k}x_k(n-1)}}{\sqrt{2\pi \sqrt{\alpha_k} x_k(n-1)}} \varphi_k(x_k(n-1)) \cdot \frac{e^{\sqrt{\alpha_k}x_k(n+1)}}{\sqrt{2\pi \sqrt{\alpha_k} x_k(n+1)}} \varphi_k(x_k(n+1))}
	{\frac{e^{2\sqrt{\alpha_k}x_k(n)}}{2\pi \sqrt{\alpha_k} x_k(n)} \varPhi_k(x_k(n))^2}\\[9pt]
	&=	\frac{x_k(n)}{\sqrt{x_k(n-1) x_k(n+1)}}
e^{\sqrt{\alpha_k } (x_k(n-1)+x_k(n+1)-2x_k(n))}\ell_k(n), \label{eq-H-geq-l}
\end{align}
where 
\begin{equation}\label{eq-L-x-phi}
	L_k(n)=\frac{\varPhi_k(x_k(n-1))\varPhi_k(x_k(n+1))}{\varphi_k(x_k(n))^2}, 
\end{equation}
and	
\begin{equation}\label{eq-l-x-phi}	
			\ell_k(n)=\frac{\varphi_k(x_k(n-1))\varphi_k(x_k(n+1))}{\varPhi_k(x_k(n))^2}.
\end{equation}

To prove \eqref{eq-H-Leq} and \eqref{eq-H-Geq}, we aim to estimate $e^{\sqrt{\alpha_k } (x_k(n-1)+x_k(n+1)-2x_k(n))}$, $L_k(n)$ and $\ell_k(n)$ in terms of $x_k(n)$.

Recalling the expression of $x_k(n)$, that is, 
\begin{equation*}
	x_k(n)=\frac{\pi\sqrt{24n-(2k+2)}}{6},
\end{equation*}
and noting that for $n\geq 2$,
\begin{align}\label{eq-x-n}
	x_k(n-1)=\sqrt{x_k(n)^2-\frac{2\pi^2}{3}},\quad x_k(n+1)=\sqrt{x_k(n)^2+\frac{2\pi^2}{3}}.
\end{align}
Then we have
\begin{align*}
	&x_k(n-1)=x_k(n)-\frac{\pi^2}{3x_k(n)}-\frac{\pi^4}{18x^3_k(n)}
      -\frac{\pi^6}{54x_k^5(n)}-\frac{5\pi^8}{648x_k^7(n)}+o\left(\frac{1}{x_k^9(n)}\right),\\[9pt]	
    &x_k(n+1)=x_k(n)+\frac{\pi^2}{3x_k(n)}-\frac{\pi^4}{18x^3_k(n)}
      +\frac{\pi^6}{54x_k^5(n)}-\frac{5\pi^8}{648x_k^7(n)}+o\left(\frac{1}{x_k^9(n)}\right).
\end{align*}
It can be readily checked that for $x_k(n)\geq 6$,
\begin{equation}\label{xkn-1}
	\tilde{w}_k(n)<x_k(n-1)<\hat{w}_k(n),
\end{equation}
and
\begin{equation}\label{xkn+1}
		\tilde{y}_k(n)<x_k(n+1)<\hat{y}_k(n), 
\end{equation}
where
\begin{align}\left\{
\begin{aligned}\label{wylabel}
	\tilde{w}_k(n)&=x_k(n)-\frac{\pi^2}{3x_k(n)}-\frac{\pi^4}{18x^3_k(n)}-\frac{\pi^6}{54x_k^5(n)}-\frac{5\pi^8}{324x_k^7(n)},\\[9pt]
	\hat{w}_k(n)&=x_k(n)-\frac{\pi^2}{3x_k(n)}-\frac{\pi^4}{18x^3_k(n)}-\frac{\pi^6}{54x_k^5(n)},\\[9pt]
\tilde{y}_k(n)&=x_k(n)+\frac{\pi^2}{3x_k(n)}-\frac{\pi^4}{18x^3_k(n)}+
\frac{\pi^6}{54x_k^5(n)}-\frac{5\pi^8}{324x_k^7(n)},\\[9pt]
	\hat{y}_k(n)&=x_k(n)+\frac{\pi^2}{3x_k(n)}-\frac{\pi^4}{18x^3_k(n)}+\frac{\pi^6}{54x_k^5(n)}.
\end{aligned}\right.
\end{align}

With these bounds of $x_k(n-1)$ and $x_k(n+1)$ in \eqref{xkn-1} and \eqref{xkn+1}, we are now in a position to bound $e^{\sqrt{\alpha_k } (x_k(n-1)+x_k(n+1)-2x_k(n))}$, $L_k(n)$ and $\ell_k(n)$.
 
First, we consider the exponential function $e^{\sqrt{\alpha_k } (x_k(n-1)+x_k(n+1)-2x_k(n))}$.
Applying \eqref{xkn-1}, \eqref{xkn+1}, and \eqref{wylabel} to the exponent of this exponential function, we find that for $x_k(n)\geq 6$,
\begin{align*}
	-\frac{\sqrt{\alpha_k}\pi^4}{9x_k(n)^3}-\frac{5\sqrt{\alpha_k}\pi^8}{162x_k(n)^7}<
\sqrt{\alpha_k}(x_k(n-1)+x_k(n+1)-2x_k(n))<-\frac{\sqrt{{\alpha_k}}\pi^4}{9x_k(n)^3},
\end{align*}
which leads us to
\begin{align}\label{eq-e-R-L}
e^{-\frac{\sqrt{\alpha_k}\pi^4}{9x_k(n)^3}-\frac{5\sqrt{\alpha_k}\pi^8}{162x_k(n)^7}}	< e^{\sqrt{\alpha_k } (x_k(n-1)+x_k(n+1)-2x_k(n))}
	< e^{-\frac{\sqrt{\alpha_k}\pi^4}{9x_k(n)^3}}.
\end{align}
We note that $s<0$, 
\begin{align}\label{eq-e-s}
	1+s<e^s<1+s+s^2.
\end{align}
Thus, by \eqref{eq-e-s}, we get that 
\begin{equation}
	\label{eq-e-geq}	
e^{-\frac{\sqrt{\alpha_k}\pi^4}{9x_k(n)^3}} < 1-\frac{\sqrt{\alpha_k}\pi^4}{9x_k(n)^3}+\frac{\alpha_k\pi^8}{81x_k(n)^6},
\end{equation}
and
\begin{equation}\label{eq-e-leq}	
e^{-\frac{\sqrt{\alpha_k}\pi^4}{9x_k(n)^3}-\frac{5\sqrt{\alpha_k}\pi^8}{162x_k(n)^7}} 	> 1-\frac{\sqrt{\alpha_k}\pi^4}{9x_k(n)^3}-\frac{5\sqrt{\alpha_k}\pi^8}{162x_k(n)^7}.
\end{equation}
Combing \eqref{eq-e-R-L}, \eqref{eq-e-geq} and  \eqref{eq-e-leq} yields that for $x_k(n)\geq 6$,
	\begin{equation}\label{eq-E-left}
		e^{\sqrt{\alpha_k } (x_k(n-1)+x_k(n+1)-2x_k(n))} > 1-\frac{\sqrt{\alpha_k}\pi^4}{9x_k(n)^3}-\frac{5\sqrt{\alpha_k}\pi^8}{162x_k(n)^7},
	\end{equation}
	and 
	\begin{equation}\label{eq-E-right} 
	e^{\sqrt{\alpha_k } (x_k(n-1)+x_k(n+1)-2x_k(n))} 
		< 1-\frac{\sqrt{\alpha_k}\pi^4}{9x_k(n)^3}+\frac{\alpha_k\pi^8}{81x_k(n)^6}.
	\end{equation}

Now we turn to estimate $L_k(n)$ and $\ell_k(n)$.
Recalling the definitions of $L_k(n)$ and $\ell_k(n)$(see \eqref{eq-L-x-phi} and \eqref{eq-l-x-phi}), we can rewrite $L_k(n)$ and $\ell_k(n)$ as follows: 
\begin{align*}
	L_k(n)=\frac{P_k(n)}{Q_k(n)},\quad  \ell_k(n)=\frac{p_k(n)}{q_k(n)},
\end{align*}
where
\begin{align*}
	P_k(n)&=x_k(n)^{14}\Big(x_k(n-1)^6-\gamma _1 x_k(n-1)^5+\gamma _2 x_k(n-1)^4+\gamma _3 x_k(n-1)^3\\
	&\qquad +\gamma _4 x_k(n-1)^2+\gamma _5 x_k(n-1)+\gamma_6\Big) 
	\Big(x_k(n+1)^6-\gamma _1 x_k(n+1)^5\\
	&\qquad +\gamma _2 x_k(n+1)^4+\gamma _3 x_k(n+1)^3+\gamma _4 x_k(n+1)^2+\gamma _5 x_k(n+1)+\gamma _6\Big), \\[9pt]
	Q_k(n)&=x_k(n)^2x_k(n-1)^6 x_k(n+1)^6
	\Big(x_k(n)^6-\gamma_1 x_k(n)^5+\gamma_2 x_k(n)^4\\
	&\qquad+\gamma_3 x_k(n)^3+\gamma_4 x_k(n)^2+\gamma_5 x_k(n)-\gamma_6 \Big)^2,\\[9pt]
	p_k(n)&=x_k(n)^{14}\Big(x_k(n-1)^6-\gamma _1 x_k(n-1)^5+\gamma _2 x_k(n-1)^4+\gamma _3 x_k(n-1)^3\\
	&\qquad +\gamma _4 x_k(n-1)^2+\gamma _5 x_k(n-1)-\gamma_6\Big) 
	\Big(x_k(n+1)^6-\gamma _1 x_k(n+1)^5\\
	&\qquad +\gamma _2 x_k(n+1)^4+\gamma _3 x_k(n+1)^3+\gamma _4 x_k(n+1)^2+\gamma _5 x_k(n+1)-\gamma _6\Big),
\end{align*}
and
\begin{align*}
	q_k(n)&=x_k(n)^2x_k(n-1)^6 x_k(n+1)^6
	\Big(x_k(n)^6-\gamma_1 x_k(n)^5+\gamma_2 x_k(n)^4\\
	&\qquad+\gamma_3 x_k(n)^3+\gamma_4 x_k(n)^2+\gamma_5 x_k(n)+\gamma_6 \Big)^2.
\end{align*}

Let
\begin{align*}
	\hat{P}_k(n)&=x_k(n)^{14}\Big(x_k(n-1)^6-\gamma _1 x_k(n-1)^4\tilde{w}_k(n)+\gamma _2 x_k(n-1)^4\\
	&\qquad +\gamma _3 x_k(n-1)^2\hat{w}_k(n)+\gamma _4 x_k(n-1)^2+\gamma _5 \hat{w}_k(n)+\gamma_6\Big) 
	\\
	&\qquad \times\Big(x_k(n+1)^6 -\gamma _1 x_k(n+1)^4\tilde{y}_k(n) +\gamma _2 x_k(n+1)^4\\
	&\qquad +\gamma _3 x_k(n+1)^2\hat{y}_k(n)+\gamma _4 x_k(n+1)^2+\gamma _5 \hat{y}_k(n)+\gamma _6\Big), 
\end{align*}
and
\begin{align*}  
	\tilde{p}_k(n)&=x_k(n)^{14}\Big(x_k(n-1)^6-\gamma _1 x_k(n-1)^4\hat{w}_k(n)+\gamma _2 x_k(n-1)^4\\
	&\qquad +\gamma _3 x_k(n-1)^2\tilde{w}_k(n)+\gamma _4 x_k(n-1)^2+\gamma _5 \tilde{w}_k(n)-\gamma_6\Big)\\
	&\qquad \times \Big(x_k(n+1)^6-\gamma _1 x_k(n+1)^4 \hat{y}_k(n) +\gamma _2 x_k(n+1)^4\\
	&\qquad +\gamma _3 x_k(n+1)^2\tilde{y}_k(n)+\gamma _4 x_k(n+1)^2+\gamma _5 \tilde{y}_k(n)-\gamma _6\Big),
\end{align*}

It can be readily checked that for $x_k(n)\geq 6$, 
\begin{align*}
	0<P_k(n)\leq \hat{P}_k(n),
\end{align*}
and
\begin{align*}
	0<  \tilde{p}_k(n) \leq p_k(n).
\end{align*}
Noting that for $x_k(n)>0$, both $Q_k(n)$ and $q_k(n)$ are positive.
Thus we have that for $x_k(n)\geq 6$,
\begin{align}\label{eq-L-P-leq}
	L_k(n)=\frac{P_k(n)}{Q_k(n)}\leq \frac{\hat{P}_k(n)}{Q_k(n)},
\end{align}
and
\begin{align}\label{eq-l-p-geq}
\ell_k(n)=\frac{p_k(n)}{q_k(n)}\geq \frac{\tilde{p}_k(n)}{q_k(n)}.
\end{align}
To bound $L_{k}(n)$ and $\ell_k(n)$ in terms of $x_k(n)$,  we shall show that for $x_k(n)\geq 110$, 
\begin{align}
	\nonumber \frac{\hat{P}_k(n)}{Q_k(n)}&\leq  1-\frac{5 \pi ^4}{8 \sqrt{\alpha_k} x_k(n)^5}+\frac{292}{\alpha_k^3 x_k(n)^6}\\[9pt]
	&=\frac{2336-5 \pi ^4 \alpha_k^{5/2} x_k(n)+8 \alpha_k^3 x_k(n)^6}{8 \alpha_k^3 x_k(n)^6},\label{eq-P-L-k-x}
\end{align}
and 
\begin{align}
	\nonumber \frac{\tilde{p}_k(n)}{q_k(n)} &\geq 1-\frac{5 \pi ^4}{8 \sqrt{\alpha_k} x_k(n)^5}-\frac{5 \pi ^4}{6 \alpha_k x_k(n)^6}-\frac{300}{\alpha_k^3 x_k(n)^6}\\[9pt]
	&=\frac{-7200-20 \pi ^4 \alpha_k^2-15 \pi ^4 \alpha_k^{5/2} x_k(n)+24 \alpha_k^3 x_k(n)^6}{24 \alpha_k^3 x_k(n)^6}.\label{eq-p-l-k-x}
\end{align}
In fact, we will prove that for $x_k(n)\geq 110$,
\begin{align}
	\nonumber &\Big(2336-5 \pi ^4 \alpha_k^{5/2} x_k(n)+8 \alpha_k^3 x_k(n)^6\Big)Q_k(n)\\[3pt]
	&\qquad \quad -8 \alpha_k^3 x_k(n)^6\hat{P}_k(n)\geq 0, \label{eq-Geq}
\end{align}
and  
\begin{align}
	\nonumber &24 \alpha_k^3 x_k(n)^6 \tilde{p}_k(n)\\[3pt]
	&\quad -q_k(n)\Big(-7200-20 \pi ^4 \alpha_k^2-15 \pi ^4 \alpha_k^{5/2} x_k(n)+24 \alpha_k^3 x_k(n)^6\Big)\geq 0. \label{eq-Leq}
\end{align}
Recalling the definitions of $Q_k(n),~\hat{P}_k(n),~\tilde{p}_k(n),~q_k(n),~ \tilde{w}_k(n),~\hat{w}_k(n),~\tilde{y}_k(n)$ and $\hat{y}_k(n)$ and 
substituting $x_k(n-1)$ and $x_k(n+1)$ with the expressions in \eqref{eq-x-n},
we find that \eqref{eq-Geq} can be rewritten as 
\begin{align}
	\sum_{j=0}^{26}a_j(k) x_k(n)^j \geq 0, \label{eq-a-L}
\end{align}
and  \eqref{eq-Leq} can be expressed as
\begin{align}
	\sum_{j=0}^{26}b_j(k) x_k(n)^j \geq 0, \label{eq-b-l}
\end{align}
where $a_j(k)$ and $b_j(k)$ are determined by $k$. Here we just list the values of $a_{26}(k)$, $a_{25}(k)$, $a_{24}(k)$, $b_{26}(k)$, $b_{25}(k)$ and $b_{24}(k)$:
\begin{align*}
	a_{26}(k)&=\frac{20}{3} \pi ^4 \alpha_k^2,\\
	a_{25}(k)&=-\frac{325  }{16}\pi ^4\alpha_k^{3/2}-\frac{4380}{\sqrt{\alpha_k}},\\
	a_{24}(k)&=\frac{595}{64} \pi ^4 \alpha_k+\frac{40515}{4 \alpha_k},\\
	b_{26}(k)&=192,\\
	b_{25}(k)&=-\frac {225} {16}\pi^4\alpha_k^{3/
      2}  - \frac {13860} {\sqrt {\alpha_k}},\\
	b_{24}(k)&=\frac{4815 }{64}\pi ^4 \alpha_k+\frac{125505}{4 \alpha_k}.
\end{align*}
It can be readily checked that for any $0\leq j \leq 23$,
\begin{align*}
	-|a_j(k)|t^j\geq -|a_{24}(k)|t^{24}
\end{align*}
holds for $t\geq 18$ and for any $0\leq j \leq 23$,
\begin{align*}
	-|b_j(k)| t^j \geq -|b_{24}(k)|t^{24},
\end{align*}
holds for $t\geq 16$.
Thus for $x_k(n)\geq 18$, we have
\begin{align}
	\nonumber \sum_{j=0}^{26}a_j(k) x_k(n)^j&\geq -\sum_{j=0}^{24}|a_{j}(k)|x_{k}(n)^j+a_{25}(k)x_k(n)^{25}+a_{26}(k) x_k(n)^{26}\\
	&\geq -25|a_{24}(k)| x_{k}(n)^{24}+a_{25}(k)x_k(n)^{25}+a_{26}(k) x_k(n)^{26}, \label{eq-Max-L}
\end{align}
and 
\begin{align}
	\nonumber \sum_{j=0}^{27}b_j(k) x_k(n)^j&\geq -\sum_{j=0}^{24}|b_{j}(k)|x_{k}(n)^j+b_{25}(k)x_k(n)^{25}+b_{26}(k) x_k(n)^{26}\\
	&\geq -25|b_{24}(k)| x_{k}(n)^{24}+b_{25}(k)x_k(n)^{25}+b_{26}(k) x_k(n)^{26}. \label{eq-Min-l}
\end{align}
For $x_k(n)\geq 110$, one can easily check that 
\begin{align}
	-25|a_{24}(k)| x_{k}(n)^{24}+a_{25}(k)x_k(n)^{25}+a_{26}(k) x_k(n)^{26}\geq 0, \label{eq-2-Max}
\end{align}
and 
\begin{align}
-25|b_{24}(k)| x_{k}(n)^{24}+b_{25}(k)x_k(n)^{25}+b_{26}(k) x_k(n)^{26}\geq 0. \label{eq-2-Min}	
\end{align}

Hence \eqref{eq-a-L} and \eqref{eq-b-l} hold for $x_k(n)\geq 110$.
So do \eqref{eq-P-L-k-x} and \eqref{eq-p-l-k-x}.
 Combining \eqref{eq-L-P-leq} and \eqref{eq-P-L-k-x}, we obtain that for $x_k(n)\geq 110$, 
\begin{equation} \label{eq-L-x}
	L_k(n)\leq 1-\frac{5 \pi ^4}{8 \sqrt{\alpha_k} x_k(n)^5}+\frac{292}{\alpha_k^3 x_k(n)^6},
\end{equation}
In view of  \eqref{eq-l-p-geq} and \eqref{eq-p-l-k-x}, we get that for $x_k(n)\geq 110$,
\begin{equation}\label{eq-l-x}
\ell_k(n) \geq 1-\frac{5 \pi ^4}{8 \sqrt{\alpha_k} x_k(n)^5}-\frac{5 \pi ^4}{6 \alpha_k x_k(n)^6}-\frac{300}{\alpha_k^3 x_k(n)^6}. 
\end{equation}

Combing \eqref{eq-H-leq-L}, \eqref{eq-E-right} and \eqref{eq-L-x}, we are lead to \eqref{eq-H-Leq}. And applying \eqref{eq-E-left} and \eqref{eq-l-x} to \eqref{eq-H-geq-l} gives us \eqref{eq-H-Geq}. This completes the proof. 
\qed

We are now ready  to prove Theorem \ref{the-Lambda-bounding}.

{\it Proof of Theorem \ref{the-Lambda-bounding}.} 
 Recall that 
\[\Lambda_k(n)=\frac{M_k(n-1)M_k(n+1)}{M_k(n)^2},\]
	as defined in \eqref{def-lambda}, where $M_k(n)$ is defined as in \eqref{def-Mk}, that is, 
\[M_k(n)=\frac{\alpha_k \pi^3}{18  {x^2_k(n)}}  I_2\left(\sqrt{\alpha_k } {x_k(n)}\right).\]
Let $H_k(n)$ be defined as in \eqref{eq-H-k}. Then we find that 
\begin{equation}
	\Lambda_k(n)=\frac{x_k(n)^4}{x_k(n-1)^2x_k(n+1)^2}H_k(n).
\end{equation}
By Lemma \ref{lem-H-bounding}, we see that for $n\geq 143296$, 

\begin{align}
	\nonumber	\Lambda_k(n)&< \frac{x_k(n)^5}{\sqrt{x_k(n-1)^5 x_k(n+1)^5}}
						\left(1- \frac{\sqrt{\alpha_k}\pi^4}{9x_k(n)^3}+\frac{\alpha_k\pi^8}{81x_k(n)^6}\right)\\[6pt]
    \label{eq-Lambda-Leq}  &\qquad \times 
    					\left( 1-\frac{5 \pi ^4}{8 			  \sqrt{\alpha_k} x_k(n)^5}+\frac{292}{\alpha_k^3 x_k(n)^6}\right),\\[9pt]
	\nonumber	\Lambda_k(n)&> \frac{x_k(n)^5}{\sqrt{x_k(n-1)^5 x_k(n+1)^5}}
						\left(1-\frac{\sqrt{\alpha_k}\pi^4}{9x_k(n)^3}-\frac{5\sqrt{\alpha_k}\pi^8}{162x_k(n)^7}\right)
						\\[6pt]
	\label{eq-Lambda-Geq}  &\qquad \times\left( 1-\frac{5 \pi ^4}{8 \sqrt{\alpha_k} x_k(n)^5}-\frac{5 \pi ^4}{6 \alpha_k x_k(n)^6}-\frac{300}{\alpha_k^3 x_k(n)^6}\right).
 \end{align}
 
 To prove \eqref{eq-Lambda-leq} and \eqref{eq-Lambda-geq},
  we first show that for $x_k(n)\geq 5$,
\begin{align}\label{eq-x-3}
	1+\frac{5 \pi ^4}{9 x(n)^4}+\frac{5\pi ^8}{18 x(n)^8}\leq \frac{x_k(n)^5}{\sqrt{x_k(n-1)^5 x_k(n+1)^5}}\leq 1+\frac{5 \pi ^4}{9 x(n)^4}+\frac{\pi ^8}{3 x(n)^8},
\end{align}
which is equivalent to 
\begin{equation}\label{eq-x-n-3}
   \left\{
	\begin{aligned}
		&x_k(n)^{20}-x_k(n-1)^{10}x_k(n+1)^{10}\left(1+\frac{5\pi ^4}{9 x_k(n)^4}+\frac{5 \pi ^8}{18 x_k(n)^8}\right)^4 \geq 0,\\[5pt]
		&x_k(n)^{20}-x_k(n-1)^{10}x_k(n+1)^{10}\left(1+\frac{5\pi ^4}{9 x_k(n)^4}+\frac{ \pi ^8}{3 x_k(n)^8}\right)^4 \leq 0.
	\end{aligned}
	\right.
\end{equation}
Recall that  
\begin{align*}
	x_k(n-1)=\sqrt{x_k(n)^2-\frac{2\pi^2}{3}},\quad x_k(n+1)=\sqrt{x_k(n)^2+\frac{2\pi^2}{3}}.
\end{align*}
Then it can be calculated that 
\begin{align}
\nonumber &x_k(n)^{20}-x_k(n-1)^{10}x_k(n+1)^{10}\left(1+\frac{5\pi ^4}{9 x_k(n)^4}+\frac{5 \pi ^8}{18 x_k(n)^8}\right)^4\\[9pt]
\nonumber
 &\quad=\frac{\pi ^{12}}{6198727824 x_k(n)^{32}}  \Big(3316191840 x_k(n)^{40}-276349320 \pi ^4 x_k(n)^{36} -36846576 \pi ^8 x_k(n)^{32}\\[9pt]
\nonumber
 & \qquad \quad 
  -673246080 \pi ^{12}x_k(n)^{28}+108795960 \pi ^{16} x_k(n)^{24} +9555975 \pi ^{20} x_k(n)^{20}\\[9pt]
\nonumber
 &\qquad \quad +61090900 \pi ^{24} x_k(n)^{16}-14404000 \pi ^{28} x_k(n)^{12}-624000 \pi ^{32} x_k(n)^8\\[9pt]
 &\qquad \quad -2080000 \pi ^{36} x_k(n)^4+640000 \pi ^{40}\Big), \label{eq-x-3-1}
\end{align}
and 
\begin{align}
	\nonumber
	&x_k(n)^{20}-x_k(n-1)^{10}x_k(n+1)^{10}\left(1+\frac{5\pi ^4}{9 x_k(n)^4}+\frac{ \pi ^8}{ x_k(n)^8}\right)^4\\[9pt]
	\nonumber
	&\quad=-\frac{\pi ^8 }{387420489 x_k(n)^{32}} \Big(86093442 x_k(n)^{44}-255091680 \pi ^4 x_k(n)^{40}+27103491 \pi ^8 x_k(n)^{36}\\[9pt]
	\nonumber
	&\qquad \quad -39917124 \pi ^{12} x_k(n)^{32}+67127778 \pi ^{16} x_k(n)^{28}-11626092 \pi ^{20} x_k(n)^{24}\\[9pt]
	\nonumber
	&\qquad \quad +6422409 \pi ^{24} x_k(n)^{20}-8161876 \pi ^{28} x_k(n)^{16}+1691040 \pi ^{32} x_k(n)^{12}\\[9pt]
	&\qquad \quad -355968 \pi ^{36} x_k(n)^8+380160 \pi ^{40} x_k(n)^4-82944 \pi ^{44}\Big). \label{eq-x-3-2}
\end{align}

Note that for $x_k(n)\geq 3$,
\begin{equation}\label{eq-1}
\left\{
	\begin{aligned}
	3316191840 x_k(n)^{40}-276349320 \pi ^4 x_k(n)^{36}-36846576 \pi ^8 x_k(n)^{32}-673246080 \pi ^{12} x_k(n)^{28}\geq  &0,\\[6pt]
	61090900 \pi ^{24} x_k(n)^{16}-14404000 \pi ^{28} x_k(n)^{12}-624000 \pi ^{32} x_k(n)^8-2080000 \pi ^{36} x_k(n)^4\geq &0,
\end{aligned}
\right.
\end{equation}
and for $x_k(n)\geq 5$,
\begin{equation}\label{eq-2}
\left\{
	\begin{aligned}
86093442 x_k(n)^{44}-255091680 \pi ^4 x_k(n)^{40}&\geq 0,\\[6pt]
27103491 \pi ^8 x_k(n)^{36} -39917124 \pi ^{12} x_k(n)^{32}&\geq 0,\\[6pt]
67127778 \pi ^{16} x_k(n)^{28}-11626092 \pi ^{20} x_k(n)^{24}&\geq 0,\\[6pt]
6422409 \pi ^{24} x_k(n)^{20}-8161876 \pi ^{28} x_k(n)^{16}&\geq 0,\\[6pt]
1691040 \pi ^{32} x_k(n)^{12}-355968 \pi ^{36} x_k(n)^8&\geq 0,\\[6pt]
380160 \pi ^{40} x_k(n)^4-82944 \pi ^{44}&\geq 0.
\end{aligned}
\right.
\end{equation}
Applying \eqref{eq-1} to \eqref{eq-x-3-1} and applying \eqref{eq-2} to \eqref{eq-x-3-2}, we find that \eqref{eq-x-n-3} holds for $x_k(n)\geq 5$, which implies \eqref{eq-x-3} is true for $x_k(n)\geq 5$. 

Combing  \eqref{eq-Lambda-Leq}, \eqref{eq-Lambda-Geq} and \eqref{eq-x-3},
we obtain that both \eqref{eq-Lambda-leq} and \eqref{eq-Lambda-geq} hold for $n\geq 143296$. The case for $2 \leq n \leq 143295$ can be directly verified by computer, and hence the proof is complete.
\qed

\section{PROOF OF THEOREM \ref{the-Theta-bounding}}

In this section, we present a proof of Theorem \ref{the-Theta-bounding}, which is based on the inequalities in the previous sections.
Let $\Theta_k(n)$ be defined as \eqref{def-Theta}, that is, 
	\begin{equation*}
			\Theta_k(n)= \frac{\Delta_k(n-1)\Delta_k(n+1)}{\Delta_k(n)^2}.
	\end{equation*}
The theorem states that for $k=1$ or $2$ and  for $x_k(n)\geq 315$, 
		\begin{align}\label{ineq-Theta-geq-1}
		\Theta_k(n)> 
		 1-\frac{\pi ^4 \sqrt{\alpha_k}}{9 x_k(n)^3}+\frac{5\pi ^4}{9 x_k(n)^4}-\frac{5 \pi ^4}{8 x_k(n)^5 \sqrt{\alpha_k}}+\frac{-\frac{300}{\alpha_k^3}-10-\frac{5 \pi ^4}{6 \alpha_k}}{x_k(n)^6},
		\end{align}
		and 
		\begin{align}\label{ineq-Theta-leq-1}
			\Theta_k(n)<
		1-\frac{\pi ^4 \sqrt{\alpha_k}}{9 x_k(n)^3}+\frac{5\pi ^4}{9 x_k(n)^4}-\frac{5 \pi ^4}{8 x_k(n)^5 \sqrt{\alpha_k}}+\frac{\frac{ \pi ^8 \alpha_k}{81}+\frac{292}{\alpha_k^3}+5}{x_k(n)^6}.
		\end{align}
{\it Proof of Theorem \ref{the-Theta-bounding}.} 
Let $g_k(n)$ be defined as in \eqref{def-g-k-n} and $G_k(n)$ be defined as in \eqref{def-G-k-n}, that is,
\begin{equation*}
\left\{
\begin{aligned}
	g_k(n)&=\frac{\left(1-\frac{1}{x_k(n-1)^6}\right)\left(1-\frac{1}{x_k(n+1)^6}\right)}{\left(1+\frac{1}{x_k(n)^6}\right)^2},\\[9pt]
		G_k(n)&=\frac{\left(1+\frac{1}{x_k(n-1)^6}\right)\left(1+\frac{1}{x_k(n+1)^6}\right)}{\left(1-\frac{1}{x_k(n)^6}\right)^2}.
\end{aligned}
\right.
\end{equation*}
		
\noindent As sated in Corollary \ref{cor-Dong-Ji-Jia}, there is a rough lower and upper bounds for $\Theta_k(n)$:
	\[
		\Lambda_k(n) g_k(n) \leq \Theta_k(n) \leq \Lambda_k(n) G_k(n).
	\]
Applying the upper and lower bounds for $\Lambda_k(n)$ in Theorem \ref{the-Lambda-bounding}, we find that for $x_k(n)\geq 152$,
\begin{align}
	\nonumber \Theta_k(n) &<
	\left(1+\frac{5\pi ^4}{9 x_k(n)^4}+\frac{ \pi ^8}{3 x_k(n)^8}\right)\left(1-\frac{\sqrt{\alpha_k}\pi^4}{9x_k(n)^3}+\frac{\alpha_k\pi^8}{81x_k(n)^6}\right)\\
	\label{eq-Theta-leq} &\qquad \times \left( 1-\frac{5 \pi ^4}{8 \sqrt{\alpha_k} x_k(n)^5}+\frac{292}{\alpha_k^3 x_k(n)^6}\right) G_k(n),
\end{align}
and 
\begin{align}
\nonumber \Theta_k(n) &> \left(1+\frac{5\pi ^4}{9 x_k(n)^4}+\frac{5 \pi ^8}{18 x_k(n)^8}\right)\left(1-\frac{\sqrt{\alpha_k}\pi^4}{9x_k(n)^3}-\frac{5\sqrt{\alpha_k}\pi^8}{162x_k(n)^7}\right)\\
	\label{eq-Theta-geq}	&\qquad \times\left( 1-\frac{5 \pi ^4}{8 \sqrt{\alpha_k} x_k(n)^5}-\frac{5 \pi ^4}{6 \alpha_k x_k(n)^6}-\frac{300}{\alpha_k^3 x_k(n)^6}\right) g_k(n).	
\end{align}

To verify \eqref{ineq-Theta-geq-1} and \eqref{ineq-Theta-leq-1}, 
we first claim that  for $x_k(n)\geq 6$, 
\begin{align}\label{ineq-g-x-6}
g_k(n) \geq 1-\frac{5}{x_k(n)^6},
\end{align}
and 
\begin{align}\label{ineq-G-x-6}
G_k(n) \leq 1+\frac{5}{x_k(n)^6}.	
\end{align}
Recall that 
\[x_k(n-1)=\sqrt{x_k(n)^2-\frac{2\pi^2}{3}},\quad x_k(n+1)=\sqrt{x_k(n)^2+\frac{2\pi^2}{3}}.
\]
Thus we obtain that 
\begin{align*}
	g_k(n)&=\frac{\left(1-\frac{1}{x_k(n-1)^6}\right)\left(1-\frac{1}{x_k(n+1)^6}\right)}{\left(1+\frac{1}{x_k(n)^6}\right)^2}\\[9pt]
	&=\frac{x_k(n)^{12} \left(\left(x_k(n)^2+\frac{2 \pi^2}{3}\right)^3-1\right) \left(\left(x_k(n)^2-\frac{2 \pi ^2}{3}\right)^3-1\right)}{\left(x_k(n)^6+1\right)^2 \left(x_k(n)^4-\frac{4 \pi ^4}{9}\right)^3},
\end{align*}
and 
\begin{align*}
	G_k(n)&=\frac{\left(1+\frac{1}{x_k(n-1)^6}\right)\left(1+\frac{1}{x_k(n+1)^6}\right)}{\left(1-\frac{1}{x_k(n)^6}\right)^2}\\[9pt]
	&=\frac{x_k(n)^{12} \left(\left(x_k(n)^2-\frac{2\pi^2} {3}\right)^3+1\right) \left(\left(x_k(n)^2+\frac{2\pi^2} {3} \right)^3+1\right)}{\left(x_k(n)^6-1\right)^2 \left(x_k(n)^4-\frac{4\pi^4}{9}\right)^3}.
\end{align*}
Hence,  we have that 
\begin{align}\label{eq-g-x-6}
	g_k(n)-\left(1-\frac{5}{x_k(n)^6}\right)
	=\frac{\phi(x_k(n))}{x_k(n)^6 \left(9 x_k(n)^4-4 \pi ^4\right)^3 \left(x_k(n)^6+1\right)^2},
\end{align}
and 
\begin{align}\label{eq-G-x-6}
	G_k(n)-\left(1+\frac{5}{x_k(n)^6}\right)=\frac{-\psi(x_k(n))}{x_k(n)^6 \left(9 x_k(n)^4-4 \pi ^4\right)^3 \left(x_k(n)^6-1\right)^2},
\end{align}
where 
\begin{align*}
	\phi(t)&=729 t^{24}-4860 \pi ^4 t^{20}+7290 t^{18}+1296 \pi ^8 t^{16}-8748 \pi ^4 t^{14}-192 \pi ^{12} t^{12}\\
	&\qquad +3645 t^{12}+3888 \pi ^8 t^{10}-4860 \pi ^4 t^8-576 \pi ^{12} t^6+2160 \pi ^8 t^4-320 \pi ^{12},
\end{align*}
and
\begin{align*}
	\psi(t)&=729 t^{24}-4860 \pi ^4 t^{20}-7290 t^{18}+1296 \pi ^8 t^{16}+8748 \pi ^4 t^{14}-192 \pi ^{12} t^{12}\\
	&\qquad +3645 t^{12}-3888 \pi ^8 t^{10}-4860 \pi ^4 t^8+576 \pi ^{12} t^6+2160 \pi ^8 t^4-320 \pi ^{12}.
\end{align*}
It can be readily checked that for $t\geq 6$,
\begin{equation*}
\left\{
\begin{aligned}
	729 t^{24}-4860 \pi ^4 t^{20}-7290 t^{18} &\geq 0,\\
	1296 \pi ^8 t^{16}+8748 \pi ^4 t^{14}-192 \pi ^{12} t^{12}+3645 t^{12}-3888 \pi ^8 t^{12}-4860 \pi ^4 t^{12} &\geq 0,\\
	576 \pi ^{12} t^6+2160 \pi ^8 t^4-320 \pi ^{12}& \geq 0,
\end{aligned}
\right.
\end{equation*}
which implies that for $t \geq 6$,
\begin{align*}
	\psi(t)\geq 0.
\end{align*}
We note that for $t\geq 3.3$, 
\begin{equation}
	\phi(t)\geq \psi(t).
\end{equation}
Thus we have that for $x_k(n)\geq 6$,
\begin{align}\label{ineq-phi-psi-0}
	\phi(x_k(n))\geq \psi(x_k(n))\geq 0.
\end{align}
Then our claim is confirmed by applying \eqref{ineq-phi-psi-0} to  \eqref{eq-g-x-6} and \eqref{eq-G-x-6} respectively.

In view of   \eqref{eq-Theta-leq}, \eqref{eq-Theta-geq}, \eqref{ineq-g-x-6} and \eqref{ineq-G-x-6}, we get that for $x_k(n)\geq 152$,
\begin{align*}
	 \Theta_k(n) &< 
	\left(1+\frac{5\pi ^4}{9 x_k(n)^4}+\frac{\pi ^8}{3 x_k(n)^8}\right)\left(1-\frac{\sqrt{\alpha_k}\pi^4}{9x_k(n)^3}+\frac{\alpha_k\pi^8}{81x_k(n)^6}\right)\\
      &\qquad \times \left( 1-\frac{5 \pi ^4}{8 \sqrt{\alpha_k} x_k(n)^5}+\frac{292}{\alpha_k^3 x_k(n)^6}\right) \left(1+\frac{5}{x_k(n)^6}\right),
\end{align*}
and 
\begin{align*}
 \Theta_k(n) &> \left(1+\frac{5\pi ^4}{9 x_k(n)^4}+\frac{5 \pi ^8}{18 x_k(n)^8}\right)\left(1-\frac{\sqrt{\alpha_k}\pi^4}{9x_k(n)^3}-\frac{5\sqrt{\alpha_k}\pi^8}{162x_k(n)^7}\right)\\
	&\qquad \times\left( 1-\frac{5 \pi ^4}{8 \sqrt{\alpha_k} x_k(n)^5}-\frac{5 \pi ^4}{6 \alpha_k x_k(n)^6}-\frac{300}{\alpha_k^3 x_k(n)^6}\right) \left(1-\frac{5}{x_k(n)^6}\right).	
\end{align*}

Thus \eqref{ineq-Theta-geq-1} and \eqref{ineq-Theta-leq-1} are justified if we can prove that for $x_k(n)\geq 315$,
\begin{align}
\nonumber	&\left(1+\frac{5\pi ^4}{9 x_k(n)^4}+\frac{ \pi ^8}{3 x_k(n)^8}\right)\left(1-\frac{\sqrt{\alpha_k}\pi^4}{9x_k(n)^3}+\frac{\alpha_k\pi^8}{81x_k(n)^6}\right)\\
\nonumber        &\qquad \times \left( 1-\frac{5 \pi ^4}{8 \sqrt{\alpha_k} x_k(n)^5}+\frac{292}{\alpha_k^3 x_k(n)^6}\right) \left(1+\frac{5}{x_k(n)^6}\right)\\
\label{ineq-1-leq}    &< 1 -\frac{\pi ^4 \sqrt{\alpha_k}}{9 x_k(n)^3}
      +\frac{5\pi^4}{9 x_k(n)^4}
      -\frac{5 \pi ^4}{8 \sqrt{\alpha_k} x_k(n)^5}
      +\frac{\frac{\pi ^8 \alpha _k}{81}+\frac{292}{\alpha_k^3}+5}{x_k(n)^6},
\end{align}
and
\begin{align}
\nonumber	&\left(1+\frac{5\pi ^4}{9 x_k(n)^4}+\frac{5 \pi ^8}{18x_k(n)^8}\right)\left(1-\frac{\sqrt{\alpha_k}\pi^4}{9x_k(n)^3}-\frac{5\sqrt{\alpha_k}\pi^8}{162x_k(n)^7}\right)\\
\nonumber	&\qquad \times\left( 1-\frac{5 \pi ^4}{8 \sqrt{\alpha_k} x_k(n)^5}-\frac{5 \pi ^4}{6 \alpha_k x_k(n)^6}-\frac{300}{\alpha_k^3 x_k(n)^6}\right) \left(1-\frac{5}{x_k(n)^6}\right)\\
\label{ineq-1-geq}	&> 1 - \frac {\pi^4\sqrt {\alpha _k}} {9 x_k (n)^3}
+\frac {5\pi^4} {9 x_k (n)^4}
-\frac {5\pi^4} {8\sqrt {\alpha _k} x_k (n)^5}
+\frac {-\frac {300} {\alpha _k^3} - \frac {5\pi^4} {6\alpha _k} - 10} {x_k (n)^6}.
\end{align}

To prove \eqref{ineq-1-leq}, we will show the difference between the left-hand and  the right-hand side of this inequality is negative for $x_k(n)\geq 315$,  that is, 
\begin{align}
	\nonumber &\left(1+\frac{5\pi ^4}{9 x_k(n)^4}+\frac{ \pi ^8}{3 x_k(n)^8}\right)\left(1-\frac{\sqrt{\alpha_k}\pi^4}{9x_k(n)^3}+\frac{\alpha_k\pi^8}{81x_k(n)^6}\right)\\
    \nonumber  &\qquad \times \left( 1-\frac{5 \pi ^4}{8 \sqrt{\alpha_k} x_k(n)^5}+\frac{292}{\alpha_k^3 x_k(n)^6}\right) \left(1+\frac{5}{x_k(n)^6}\right)\\
    \label{ineq-2-leq} &-\left( 1 -\frac{\pi ^4 \sqrt{\alpha_k}}{9 x_k(n)^3}
      +\frac{5\pi^4}{9 x_k(n)^4}
      -\frac{5 \pi ^4}{8 \sqrt{\alpha_k} x_k(n)^5}
      +\frac{\frac{\pi ^8 \alpha _k}{81}+\frac{292}{\alpha_k^3}+5}{x_k(n)^6}\right)<0.
\end{align}
We note that  the left-hand side of the inequality \eqref{ineq-2-leq} can be expressed as follows:
\begin{equation*}
	-\frac{J(x_k(n))}{5832 \alpha_k^3 x_k(n)^{26}},
\end{equation*}
where
\begin{equation}
	J(x_k(n))=\sum_{j=0}^{19} c_j(k) x_k(n)^j.
\end{equation}
Here we just list the values of $c_{19}(k),c_{18}(k),c_{17}(k)$:
\begin{align*}
	c_{19}(k)&=360 \pi ^8 \alpha_k^{7/2},\\
	c_{18}(k)&=-2349 \pi ^8 \alpha_k^3,\\
	c_{17}(k)&=2025 \pi ^8 \alpha_k^{5/2}+3240 \pi ^4 \alpha_k^{7/2}+189216 \pi ^4 \sqrt{\alpha_k}.
\end{align*}

Similarly, to justify \eqref{ineq-1-geq}, we aim to prove that for $x_k(n)\geq 315$,
\begin{align}
	\nonumber	&\left(1+\frac{5\pi ^4}{9 x_k(n)^4}+\frac{ 5\pi ^8}{18 x_k(n)^8}\right)\left(1-\frac{\sqrt{\alpha_k}\pi^4}{9x_k(n)^3}-\frac{5\sqrt{\alpha_k}\pi^8}{162x_k(n)^7}\right)\\
\nonumber	&\qquad \times\left( 1-\frac{5 \pi ^4}{8 \sqrt{\alpha_k} x_k(n)^5}-\frac{5 \pi ^4}{6 \alpha_k x_k(n)^6}-\frac{300}{\alpha_k^3 x_k(n)^6}\right) \left(1-\frac{5}{x_k(n)^6}\right)\\
\label{ineq-2-geq}	&-\left( 1 - \frac {\pi^4\sqrt {\alpha _k}} {9 x_k (n)^3}
+\frac {5\pi^4} {9 x_k (n)^4}
-\frac {5\pi^4} {8\sqrt {\alpha _k} x_k (n)^5}
+\frac {-\frac {300} {\alpha _k^3} - \frac {5\pi^4} {6\alpha _k} - 10} {x_k (n)^6}\right)>0.
\end{align}
And the left-hand side of \eqref{ineq-2-geq} can be rewritten as follows:
\begin{equation*}
	\frac{K(x_k(n))}{69984\alpha_k^3 x_k (n)^{27}},
\end{equation*}
where 
\begin{equation}
	K(x_k(n))=\sum_{j=0}^{21} d_j(k) x_k(n)^j. 
\end{equation}
Here we also list the values of the  first three coefficients of $K(x_k(n))$:
\begin{align*}
	d_{21}(k)&=349920\alpha _k^3,\\
	d_{20}(k)&=-6480 \pi ^8 \alpha _k^{7/2},\\
	d_{19}(k)&=24300 \pi ^8 \alpha _k^3.
\end{align*}

Since $\alpha_k$ is positive for $k=1$ or $2$, to conclude \eqref{ineq-2-leq} and \eqref{ineq-2-geq}, we turn to prove the polynomials $J(x_k(n))$ and $K(x_k(n))$ are positive for $x_k(n)\geq 315$.

By the same method applied in estimating $L_k(n)$ and $\ell_k(n)$ (see, Sec. \ref{sec-Lambda}), we find that $J(x_k(n))$ and $K(x_k(n))$ are positive if both
\begin{equation}\label{ineq-J-geq}
	J(x_k(n)) \geq \left(-17 |c_{17}(k)|+c_{18}(k)x_k(n)+ c_{19}(k) x_{k}(n)^{2}\right)x_k(n)^{17}>0,
\end{equation}
and 
\begin{equation}\label{ineq-K-geq}
	K(x_k(n)) \geq \left(-19 |d_{19}(k)| +d_{20}(k)x_k(n)+ d_{21}(k) x_{k}(n)^2\right)x_k(n)^{19}>0.
\end{equation}
hold for $x_k(n)\geq 315$.

Note that for $x_k(n)\geq 315$,
\begin{equation}\label{ineq-c-1}
	-17 |c_{17}(k)|+c_{18}(k)x_k(n)+ c_{19}(k) x_{k}(n)^{2}>0,
\end{equation}
and 
\begin{equation}\label{ineq-d-1}
	-19 |d_{19}(k)| +d_{20}(k)x_k(n)+ d_{21}(k) x_{k}(n)^{2}>0.
\end{equation}
Moreover, for any $0\leq j \leq 16$, numerical evidence indicates that for $x_k(n)\geq 7$,
\begin{equation}\label{ineq-c-2}
	-|c_j(k)|x_k(n)^j>-|c_{17}(k)| x_k(n)^{17},
\end{equation}
and for any $0\leq j \leq 18$,
\begin{equation}\label{ineq-d-2}
	-|d_j(k)|x_k(n)^j>-|d_{19}(k)| x_k(n)^{19},
\end{equation}
holds for $x_k(n)\geq 3$. Thus we have that for $x_k(n)\geq 315$,
\begin{equation}\label{ineq-J-1}
	J(x_k(n)) \geq \left(-17 |c_{17}(k)|+c_{18}(k)x_k(n)+ c_{19}(k) x_{k}(n)^{2}\right)x_k(n)^{17},
\end{equation}
and 
\begin{equation}\label{ineq-K-1}
		K(x_k(n)) \geq \left(-19 |d_{19}(k)| +d_{20}(k)x_k(n)+ d_{21}(k) x_{k}(n)^2\right)x_k(n)^{19}.
\end{equation}
Combing \eqref{ineq-c-1} and \eqref{ineq-J-1} bears \eqref{ineq-J-geq}. This proves \eqref{ineq-2-leq}. Similarly, \eqref{ineq-K-geq} is confirmed by  joining \eqref{ineq-d-1} with \eqref{ineq-K-1}. This proofs \eqref{ineq-2-geq}.

Hence we conclude that  both \eqref{ineq-Theta-geq-1} and \eqref{ineq-Theta-leq-1} hold for $x_k(n) \geq 315$, or equivalently, for $n \geq 15081$.
This completes the proof. 
\qed

\section{Proof of Theorem \ref{the-Theta-mono} and Theorem \ref{the-Theta-k}}

In this section, we will present a proof of Theorem \ref{the-Theta-mono}. And by this theorem and the following lemma, we also give a proof of Theorem \ref{the-Theta-k}. 
	\begin{lemma}\label{lem-u-v}
		Let $u$ and $v$ be two positive real numbers such that $\frac{15}{16} \leq u <v <1$. 
		If  
			\begin{equation}\label{eq-u-v-1}
				u+\sqrt{(1-u)^{3}}>v,
			\end{equation}
		then we have 
			\begin{equation}\label{eq-u-v-2}
				4(1-u)(1-v)-(1-u v)^2 > 0.
			\end{equation}
	\end{lemma}
	
	{\it Proof:}
	Let $f(t)$ be the quadratic function defined by
	\begin{equation*}
		f(t)=4(1-u)(1-t)-(1-u t)^2,
	\end{equation*}
	which can be rewritten as 
	\begin{equation*}
		f(t)=-u^2 t^2+(6 u-4)t-4u+3.
	\end{equation*}
	Then \eqref{eq-u-v-2} says $f(v)>0$.
	
	It should be noted that the equation $f(t)=0$ has two real roots:
	\begin{equation*}
		t_1(u)=\frac{3 u-2\sqrt{(1-u)^3}-2}{u^2},
	\end{equation*}
	and 
	\begin{equation*}
		t_2(u)=\frac{3 u+2\sqrt{(1-u)^3}-2}{u^2}.
	\end{equation*}
	Thus $f(t)>0$ holds if and only if $ t_1(u)<t<t_2(u)$. 
	For $\frac{15}{16} \leq u <1$,  we find that 
	\begin{equation*}
	\left\{
		\begin{aligned}
		f(0)&=3-4u<0,\\
		f(u)&=(1-u)^3(u+3)>0,\\
		f(1)&=-(1-u)^2<0.
	  \end{aligned}
	\right.
	\end{equation*}
	which implies that  
	\begin{equation*}
		0<t_1(u)<u<t_2(u)<1.
	\end{equation*}
	Since $u<v<1$, 
	to prove \eqref{eq-u-v-2}, or equivalently, $f(v)>0$,
	it is sufficient to show that
	\begin{equation}\label{eq-v-t}
		v<t_2(u).
	\end{equation}
	This goal can be achieved by finding a lower bound for $t_2(u)$. 
	More precisely, we  show that for $\frac{15}{16}\leq u<1$,
	\begin{equation}\label{ineq-t-2}
		u+\sqrt{(1-u)^3}<t_2(u).
	\end{equation}
	Let 
	\begin{equation}\label{eq-def-tau}
		\tau(u)=t_2(u)-u-\sqrt{(1-u)^3}.
	\end{equation}
	Then \eqref{ineq-t-2} says that for $\frac{15}{16}\leq u<1$,
	\begin{equation}\label{ineq-tau}
	\tau(u)>0.
	\end{equation}
	Setting $u=1-s^2$, \eqref{eq-def-tau} becomes 
	\begin{equation}
		\tau(u)=\frac{\left(-2 s+\sqrt{5}-1\right) s^3 \left(2 s+\sqrt{5}+1\right)}{4 (s+1)^2},
	\end{equation}	
	which is clearly positive for $0<s\leq \frac{1}{4}$.
	This proves \eqref{ineq-tau}, or equivalently, \eqref{ineq-t-2}.
	In view of \eqref{eq-u-v-1} and \eqref{ineq-t-2}, we arrive at \eqref{eq-v-t}.
	The proof is completed.
	\qed
	
	We now give the proof of Theorem \ref{the-Theta-mono}.
	
	{\it Proof of Theorem \ref{the-Theta-mono}:}
	According to Theorem \ref{the-Theta-bounding}, we have that for $k=1$ or $2$ and  for $x_k(n)\geq 315$,
	\begin{align}\label{eq-theta-n+1}
		\Theta_k(n+1) > 
		 1-\frac{\pi ^4 \sqrt{\alpha_k}}{9 x_k(n+1)^3}+\frac{5\pi ^4}{9x_k(n+1)^4}-\frac{5 \pi ^4}{8 x_k(n+1)^5 \sqrt{\alpha_k}}+\frac{-\frac{300}{\alpha_k^3}-10-\frac{5 \pi ^4}{6 \alpha_k}}{x_k(n+1)^6},
	\end{align}
		and 
		\begin{align}\label{eq-theta-n-<}
			\Theta_k(n)<
		1-\frac{\pi ^4 \sqrt{\alpha_k}}{9 x_k(n)^3}+\frac{5\pi ^4}{9 x_k(n)^4}-\frac{5 \pi ^4}{8 x_k(n)^5 \sqrt{\alpha_k}}+\frac{\frac{ \pi ^8 \alpha_k}{81}+\frac{292}{\alpha_k^3}+5}{x_k(n)^6}.
		\end{align}
		
Note that for $n\geq 2$,
\begin{equation}\label{eq-x-3456}
\left\{
	\begin{aligned}
		\frac{1}{x_k(n+1)^3}&< \frac{1}{x_k(n)^3}-\frac{\pi ^2}{2x_k(n)^5},\\[9pt]
	\frac{1}{x_k(n+1)^4}&>\frac{1}{x_k(n)^4}-\frac{4 \pi ^2}{3 x_k(n)^6},\\[9pt]
	\frac{1}{x_k(n+1)^5}&<\frac{1}{x_k(n)^5},\\
	\frac{1}{x_k(n+1)^6}&<\frac{1}{x_k(n)^6}.
	\end{aligned}
	\right.
\end{equation}
Applying \eqref{eq-x-3456} to \eqref{eq-theta-n+1}, we have that for $x_k(n)\geq 315$,
	\begin{align}
	\nonumber	\Theta_k(n+1) &> 
		 1-\frac{\pi ^4 \sqrt{\alpha_k}}{9 x_k(n+1)^3}+\frac{5\pi ^4}{9 x_k(n+1)^4}-\frac{5 \pi ^4}{8 x_k(n+1)^5 \sqrt{\alpha_k}}+\frac{-\frac{300}{\alpha_k^3}-10-\frac{5 \pi ^4}{6 \alpha_k}}{x_k(n+1)^6}\\[9pt]
\nonumber		& >1-\frac{\pi ^4 \sqrt{\alpha_k}}{9} \left(\frac{1}{x_k(n)^3}-\frac{\pi ^2}{2x_k(n)^5}\right) 
+\frac{5\pi ^4}{9} \left(\frac{1}{x_k(n)^4}-\frac{4 \pi ^2}{3 x_k(n)^6}\right)\\[9pt]
\nonumber &\qquad -\frac{5 \pi ^4}{8 x_k(n)^5 \sqrt{\alpha_k}}+\frac{-\frac{300}{\alpha_k^3}-10-\frac{5 \pi ^4}{6 \alpha_k}}{x_k(n)^6}\\[9pt]
&=1
-\frac{\pi ^4 \sqrt{\alpha _k}}{9 x_k(n)^3}
+\frac{5\pi ^4}{9 x_k(n)^4}
+\frac{\frac{\pi ^6 \sqrt{\alpha _k}}{18} -\frac{5 \pi ^4}{8 \sqrt{\alpha _k}}}{x_k(n)^5}
+\frac{-\frac{300}{\alpha _k^3}-\frac{5 \pi ^4}{6 \alpha _k}-\frac{20 \pi ^6}{27}-10}{x_k(n)^6}. \label{eq-theta-n+1-leq}
	\end{align}

In view of \eqref{eq-theta-n-<}	and \eqref{eq-theta-n+1-leq}, we find that for $x_k(n)\geq 315$, 
\begin{align}
\nonumber	\Theta_k(n+1)-\Theta_k(n)&>
	\Bigg(1
-\frac{\pi ^4 \sqrt{\alpha _k}}{9 x_k(n)^3}
+\frac{5\pi ^4}{9 x_k(n)^4}
+\frac{\frac{\pi ^6 \sqrt{\alpha _k}}{18} -\frac{5 \pi ^4}{8 \sqrt{\alpha _k}}}{x_k(n)^5}
+\frac{-\frac{300}{\alpha _k^3}-\frac{5 \pi ^4}{6 \alpha _k}-\frac{20 \pi ^6}{27}-10}{x_k(n)^6}\Bigg)\\[9pt]
\nonumber &\qquad  -\Bigg(1-\frac{\pi ^4 \sqrt{\alpha_k}}{9 x_k(n)^3}
+\frac{5\pi ^4}{9 x_k(n)^4}-\frac{5 \pi ^4}{8 x_k(n)^5 \sqrt{\alpha_k}}+\frac{\frac{ \pi ^8 \alpha_k}{81}+\frac{292}{\alpha_k^3}+5}{x_k(n)^6}\Bigg)\\[9pt]
&=\frac{\pi ^6 \sqrt{\alpha _k}}{18 x_k(n)^5}+\frac{-\frac{\pi ^8 \alpha _k}{81} -\frac{5 \pi ^4}{6 \alpha _k}-\frac{592}{\alpha _k^3}-\frac{20 \pi ^6}{27}-15}{x_k(n)^6}. \label{eq-theta-dif}
\end{align}
Since for $x_k(n)\geq 14$,
\begin{equation}
	\frac{\pi ^6 \sqrt{\alpha _k}}{18 x_k(n)^5}+\frac{-\frac{\pi ^8 \alpha _k}{81} -\frac{5 \pi ^4}{6 \alpha _k}-\frac{592}{\alpha _k^3}-\frac{20 \pi ^6}{27}-15}{x_k(n)^6}>0,
\end{equation}
then we get that for $k=1$ or $2$ and  $x_k(n)\geq 315$, or equivalently, for $n\geq 15081 $,
\[
	\Theta_k(n+1)-\Theta_k(n)>0.
\]
For $k=1$ or $2$, the case $n\leq 15080$ can be checked by computer.
This completes the proof.
\qed

	We are now ready to prove Theorem \ref{the-Theta-k}.

	{\it Proof of Theorem \ref{the-Theta-k}:}  
	We begin by recalling that for the log-concavity of the broken $k$-diamond function ($k=1$, or $2$), which implies that $\Theta_k(n)<1$.
	Then by Lemma \ref{lem-u-v} and Theorem \ref{the-Theta-mono}, we need only to show that 
	\begin{equation}
		\frac{15}{16} \leq \Theta_k(n),
	\end{equation}
	and 
	\begin{equation}\label{eq-Theta-3}
		\Theta_k(n+1)<\Theta_k(n)+\sqrt{(1-\Theta_k(n))^3}.
	\end{equation}

	By Theorem \ref{the-Theta-bounding}, we see that for $x_k(n)\geq 315$,
	\begin{equation}
	\Theta_k(n)>	1 - \frac {\pi^4\sqrt {\alpha _k}} {9 x_k (n)^3}
+\frac {5\pi^4} {9x_k (n)^4}
-\frac {5\pi^4} {8\sqrt {\alpha _k} x_k (n)^5}
+\frac {-\frac {300} {\alpha _k^3} - \frac {5\pi^4} {6\alpha _k} - 10} {x_k (n)^6}.
	\end{equation}
Note that for $x_k(n)\geq 0$,
\begin{align*}
	& 1 - \frac {\pi^4\sqrt {\alpha _k}} {9 x_k (n)^3}
+\frac {5\pi^4} {9 x_k (n)^4}
-\frac {5\pi^4} {8\sqrt {\alpha _k} x_k (n)^5}
+\frac {-\frac {300} {\alpha _k^3} - \frac {5\pi^4} {6\alpha _k} - 10} {x_k (n)^6}\\
& \qquad >1 - \frac {\pi^4\sqrt {\alpha _k}} {9 x_k (n)^3}
-\frac {5\pi^4} {8\sqrt {\alpha _k} x_k (n)^3}
+\frac {-\frac {300} {\alpha _k^3} - \frac {5\pi^4} {6\alpha _k} - 10} {x_k (n)^3}.
\end{align*}
It can be easily checked that for $x_k(n)\geq 13$,
\begin{align*}
	1 - \frac {\pi^4\sqrt {\alpha _k}} {9 x_k (n)^3}
-\frac {5\pi^4} {8\sqrt {\alpha _k} x_k (n)^3}
+\frac {-\frac {300} {\alpha _k^3} - \frac {5\pi^4} {6\alpha _k} - 10} {x_k (n)^3}\geq \frac{15}{16}.
\end{align*}
Then we obtain that for $x_k(n)\geq 315$,
\begin{align*}
	\Theta_k(n)\geq \frac{15}{16}.
\end{align*}

Now we shall show \eqref{eq-Theta-3}. 
By Theorem \ref{the-Theta-bounding}, we have that for $x_k(n)\geq 315$,
\begin{align}
\nonumber	&\Theta_k(n+1)-\Theta_k(n)\\
\nonumber	&\qquad <1-\frac{\pi ^4 \sqrt{\alpha_k}}{9 x_k(n+1)^3}+\frac{5\pi ^4}{9 x_k(n+1)^4}-\frac{5 \pi ^4}{8 x_k(n+1)^5 \sqrt{\alpha_k}}
	+\frac{\frac{ \pi ^8 \alpha_k}{81}+\frac{292}{\alpha_k^3}+5}{x_k(n+1)^6}\\[9pt]
\nonumber	&\qquad \quad -\left(1-\frac{\pi ^4 \sqrt{\alpha_k}}{9 x_k(n)^3}+\frac{5\pi ^4}{9 x_k(n)^4}-\frac{5 \pi ^4}{8 x_k(n)^5 \sqrt{\alpha_k}}
	+\frac{-\frac{300}{\alpha_k^3}-10-\frac{5 \pi ^4}{6 \alpha_k}}{x_k(n)^6}
	\right)\\[9pt]
\nonumber	&\qquad =\frac{\pi ^4 \sqrt{\alpha_k}}{9 } \left(\frac{1}{x_k(n)^3}-\frac{1}{x_k(n+1)^3}\right)-\frac{5\pi ^4}{9} \left(\frac{1}{x_k(n)^4}-\frac{1}{x_k(n+1)^4}\right)\\[9pt]
\label{eq-Theta-n+1-Theta-n}	&\qquad \quad +\frac{5 \pi ^4}{8 \sqrt{\alpha_k}} \left(\frac{1}{x_k(n)^5}-\frac{1}{x_k(n+1)^5}\right)
	+\frac{\frac{ \pi ^8 \alpha_k}{81}+\frac{292}{\alpha_k^3}+5}{x_k(n+1)^6}
	+\frac{\frac{300}{\alpha_k^3}+10+\frac{5 \pi ^4}{6 \alpha_k}}{x_k(n)^6},
\end{align}

Note that for  $x_k(n)\geq 4$, 
\begin{equation}\label{ineq-3456}
	\left\{
	\begin{aligned}
		&\frac{1}{x_k(n)^3}-\frac{1}{x_k(n+1)^3}<\frac{4 \pi ^2}{x_k(n)^5},\\[9pt]
		&\frac{1}{x_k(n+1)^4}-\frac{1}{x_k(n)^4}< 0,\\[9pt]
		&\frac{1}{x_k(n)^5}-\frac{1}{x_k(n+1)^5}<\frac{1}{x_k(n)^5}.
	\end{aligned}
	\right.
\end{equation}
And for $x_k(n)\geq 38$, it can be proved that 
\begin{equation}\label{eq-x-n-5-6}
	\frac{\frac{ \pi ^8 \alpha_k}{81}+\frac{292}{\alpha_k^3}+5}{x_k(n+1)^6}
	+\frac{\frac{300}{\alpha_k^3}+10+\frac{5 \pi ^4}{6 \alpha_k}}{x_k(n)^6}
<\frac{\frac{\pi ^8 \alpha _k}{81}+\frac{5 \pi ^4}{6 \alpha _k}+\frac{592}{\alpha _k^3}+15}{x_k(n)^6}<\frac{\pi^2}{x_k(n)^5}.
\end{equation}

Applying \eqref{ineq-3456} and \eqref{eq-x-n-5-6} to \eqref{eq-Theta-n+1-Theta-n}, we get that for $x_k(n)\geq 315$,
\begin{align}
\nonumber	\Theta_k(n+1)-\Theta_k(n) & <\frac{\pi ^4 \sqrt{\alpha_k}}{9 } \times  \frac{4\pi ^2}{x_k(n)^5}+\frac{5 \pi ^4}{8 \sqrt{\alpha_k} x_k(n)^5} 
	+\frac{\pi^2}{x_k(n)^5}\\[9pt]
	&=\frac{\pi ^2 \left(72 \sqrt{\alpha _k}+32 \pi ^4 \alpha _k+45 \pi ^2\right)}{72 \sqrt{\alpha _k} x_k(n)^5}. \label{eq-Theta-n+1-Theta-n-leq}
\end{align}

Since for $x_k(n)\geq 315$, 
\begin{align*}
	1-\Theta_k(n)&> \frac{\pi ^4 \sqrt{\alpha_k}}{9 x_k(n)^3}
	-\frac{5\pi ^4}{9 x_k(n)^4}
	+\frac{5 \pi ^4}{8 x_k(n)^5 \sqrt{\alpha_k}}
	-\frac{\frac{ \pi ^8 \alpha_k}{81}+\frac{292}{\alpha_k^3}+5}{x_k(n)^6}\\[9pt]
	&>\frac{\pi ^4 \sqrt{\alpha_k}}{9 x_k(n)^3}
	-\frac{5\pi ^4}{9 x_k(n)^4}
	-\frac{\frac{ \pi ^8 \alpha_k}{81}+\frac{292}{\alpha_k^3}+5}{x_k(n)^4}
\end{align*}
and for  $x_k(n)\geq 181$, it can be checked that 
\begin{align*} 
	-\frac{5\pi ^4}{9 x_k(n)^4}
	-\frac{\frac{ \pi ^8 \alpha_k}{81}+\frac{292}{\alpha_k^3}+5}{x_k(n)^4}>-\frac{2}{x_k(n)^3},
\end{align*}
thus we obtain that for $x_k(n)\geq 315$,
\begin{align*}
	1-\Theta_k(n)>\frac{\pi ^4 \sqrt{\alpha_k}}{9x_k(n)^3}-\frac{2}{x_k(n)^3}>0.
\end{align*}
Hence we have that for $x_k(n)\geq 315$,
\begin{equation}\label{ineq-1-theta-3}
\sqrt{\left(1-\Theta_k(n)\right)^3}> \frac{\sqrt{\left(\pi ^4 \sqrt{\alpha_k}-18\right)^3}}{27x_k(n)^{9/2}}.
\end{equation}
Furthermore, it can be easily checked that for $x_k(n) \geq 161$,
\begin{equation}\label{eq-9/2-5}
	\frac{\sqrt{\left(\pi ^4 \sqrt{\alpha_k}-9\right)^3}}{27x_k(n)^{9/2}}
	>
	\frac{\pi ^2 \left(72 \sqrt{\alpha _k}+32 \pi ^4 \alpha _k+45 \pi ^2\right)}{72 \sqrt{\alpha _k} x_k(n)^5}.
\end{equation}
Combing \eqref{eq-Theta-n+1-Theta-n-leq}, \eqref{ineq-1-theta-3} and \eqref{eq-9/2-5}, we obtain that for $x_k(n)\geq 315$, or equivalently, for $n\geq 15081$,
\begin{equation}
	\Theta_k(n+1)-\Theta_k(n)<\sqrt{(1-\Theta_k(n))^3},
\end{equation}
which confirms \eqref{eq-Theta-3}. 
The case for $6 \leq n \leq 15080$ can be directly checked, and  
 this completes the proof.
\qed

\end{document}